\documentclass[12pt]{amsart}
\usepackage{a4wide}
\usepackage[latin9]{inputenc}
\usepackage{amsmath}
\usepackage{amsfonts}
\usepackage{amssymb}
\usepackage{amsthm}
\usepackage{subfigure}
\newtheorem{theo}{Theorem}[section]
\newtheorem{prop}[theo]{Proposition}
\newtheorem{coro}[theo]{Corollary}
\newtheorem{lemm}[theo]{Lemma}

\theoremstyle{definition}

\theoremstyle{remark}
\newtheorem{rema}[theo]{Remark}
\newcommand{\Op}{\operatorname{Op}}

\newcommand{\nwc}{\newcommand}
\nwc{\eps}{\epsilon}
\nwc{\ep}{\epsilon}
\nwc{\vareps}{\varepsilon}
\nwc{\Oph}{\operatorname{Op}_\hbar}
\nwc{\la}{\langle}
\nwc{\ra}{\rangle}

\nwc{\mf}{\mathbf} 
\nwc{\blds}{\boldsymbol} 
\nwc{\ml}{\mathcal} 

\nwc{\defeq}{\stackrel{\rm{def}}{=}}

\nwc{\cE}{\ml{E}}
\nwc{\cN}{\ml{N}}
\nwc{\cO}{\ml{O}}
\nwc{\cP}{\ml{P}}
\nwc{\cU}{\ml{U}}
\nwc{\cV}{\ml{V}}
\nwc{\cW}{\ml{W}}
\nwc{\tU}{\widetilde{U}}
\nwc{\IN}{\mathbb{N}}
\nwc{\IR}{\mathbb{R}}
\nwc{\IZ}{\mathbb{Z}}
\nwc{\IC}{\mathbb{C}}
\nwc{\IT}{\mathbb{T}}
\nwc{\tP}{\widetilde{P}}
\nwc{\tPi}{\widetilde{\Pi}}
\nwc{\tV}{\widetilde{V}}
\nwc{\supp}{\operatorname{supp}}
\nwc{\rest}{\restriction}

\begin{document}

\title[Perturbation of the semiclassical Schr\"odinger equation]{Perturbation of the semiclassical Schr\"odinger equation on negatively curved surfaces}

\author[Suresh Eswarathasan]{Suresh Eswarathasan}
\author[Gabriel Rivi\`ere]{Gabriel Rivi\`ere}

\address{Institut des Hautes \'Etudes Scientifiques, Le Bois-Marie 35, route de Chartres, 91440 Bures-sur-Yvette, France}
\address{Department of Mathematics and Statistics, McGill University, Montr\'eal, Canada}

\email{suresh@ihes.fr, suresh@math.mcgill.ca}

\address{Laboratoire Paul Painlev\'e (U.M.R. CNRS 8524), U.F.R. de Math\'ematiques, Universit\'e Lille 1, 59655 Villeneuve d'Ascq Cedex, France}
\email{gabriel.riviere@math.univ-lille1.fr}

\begin{abstract} 
We consider the semiclassical Schr\"odinger equation on a compact negatively curved surface. For any sequence of initial data microlocalized on the unit 
cotangent bundle, we look at the quantum evolution (below the Ehrenfest time) under small perturbations of the Schr\"odinger equation, and we prove that, in the semiclassical 
limit and for typical perturbations, the solutions become equidistributed on the unit cotangent bundle.
\end{abstract} 

\keywords{Hyperbolic dynamical systems, Semiclassical analysis, Quantum chaos, Geodesic flows, Semiclassical measures}
\subjclass[2010]{37D40, 58J51, 81Q50, 35Q41}

\maketitle

\section{Introduction}

Let $(M,g)$ be a smooth ($\ml{C}^{\infty}$), compact, connected and Riemannian manifold without boundary of dimension $d$. The purpose of this article is to study the long time dynamics of the Schr\"odinger equation on $M$
\begin{equation}\label{e:schrodinger}
\imath\hbar\frac{\partial u_{\hbar}}{\partial t}=\hat{P}_0(\hbar) u_{\hbar},\qquad u_{\hbar}\rceil_{t=0}=\psi_{\hbar}\in
L^{2}(M),
\end{equation}
where $\hbar>0$, and $\hat{P}_0(\hbar):=-\frac{\hbar^2\Delta_g}{2}$, with $\Delta_g$ the Laplace Beltrami operator induced by the Riemannian metric $g$. One knows that the semiclassical properties of this equation are related to the properties of the underlying Hamiltonian system which is in this case the geodesic flow $G_0^t$ acting on the cotangent bundle $T^*M$. In the following, we will be interested by initial data which are microlocalized near a fixed energy layer of the classical flow. Precisely, we will suppose that $(\psi_{\hbar})_{\hbar>0}$ is a sequence of normalized states in $L^2(M)$ which satisfies the following oscillation assumption:
\begin{equation}\label{e:freq-assumption-unit}
\forall\delta_0>0,\ \ 
\limsup_{\hbar\rightarrow 0}\left\|\mathbf{1}_{[1/2-\delta_0,1/2+\delta_0]}\left(\hat{P}_0(\hbar)\right)\psi_{\hbar}-\psi_{\hbar}\right\|_{L^2(M)}=0.
\end{equation}
In other words, it means that the initial data are microlocalized on the unit cotangent bundle
$$T_{1/2}^*M:=\left\{(x,\xi)\in T^*M:p_0(x,\xi):=\frac{\|\xi\|_x^2}{2}=\frac{1}{2}\right\}.$$
For a given sequence of initial data $(\psi_{\hbar})_{\hbar>0}$, there are many ways to study the properties of the solutions $(u_{\hbar})_{\hbar>0}$ of~\eqref{e:schrodinger}. We will study in this article the following measures on $M$ associated to the solutions of~\eqref{e:schrodinger}:
$$\forall a\in\ml{C}^{0}(M),\ \int_{M}a(x)d\nu_{\hbar}(t)(x):=\int_{M}a(x)|u_{\hbar}(t,x)|^2d\text{vol}_M(x),$$
where $\text{vol}_{M}$ is the measure induced by the Riemannian metric on $M$. In fact, we will study a slightly more general quantity, namely
$$\forall a\in\ml{C}^{\infty}_c(T^*M),\ \mu_{\hbar}(t)(a):=\left\la \Oph(a)u_{\hbar}(t),u_{\hbar}(t)\right\ra_{L^2(M)},$$
where $\Oph(a)$ is a pseudodifferential operator of symbol $a$~\cite{Zw12} -- see also appendix~\ref{a:pdo} for a brief reminder. This quantity defines a distribution on $T^*M$, which can be extended to test functions of the form $a(x)\in\ml{C}^{\infty}(M)$ (in this case, we recover $\nu_{\hbar}(t)$). If we consider a fixed time $t\in\IR$, a simple manifestation of the relation between the classical and the quantum dynamics  is given by the Egorov Theorem which implies, in this setting, that
\begin{equation}\label{e:dist-phase-space}\forall a\in\ml{C}^{\infty}_c(T^*M),\ \mu_{\hbar}(t)(a)=\mu_{\hbar}(0)(a\circ G_0^t)+\ml{O}_t(\hbar).\end{equation}
This kind of relation still holds, modulo the fact that the remainder is slightly bigger, up to times of order $\kappa_0|\log\hbar|$~\cite{BaGrPa99, BoRo02, Zw12}, 
where $\kappa_0>0$ is some geometric constant which is related to the Hamiltonian $p_0$ and the support of $a$. Such a time is often called the Ehrenfest time. Without any additional 
geometric assumptions, understanding the dynamics of the Schr\"odinger equation beyond this logarithmic threshold is a delicate problem.

Our aim in this article is to consider the situation where the underlying classical system, being in our case the geodesic flow, enjoys some chaotic features. 
We will not try to give information on the dynamics of the Schr\"odinger equation for times longer than the Ehrenfest time. 
Our purpose is to understand how ``small perturbations'' of the equation~\eqref{e:schrodinger} alter the behavior of the solutions associated to a fixed sequence 
of initial data $(\psi_{\hbar})_{\hbar>0}$, even for short logarithmic times. In some sense, this is related to questions arising in the physics literature concerning the problem of the quantum Loschmidt echo that we will discuss below. Before being more precise on the questions that we want to address, we will first recall some results on the dynamics of the unperturbed Schr\"odinger equation when the corresponding classical system has chaotic properties such as ergodicity.


\subsection{The case of the unperturbed Schr\"odinger equation}\label{ss:unperturbed}

A natural case to consider is the situation where we make the additional assumption that $(\psi_{\hbar})_{\hbar>0}$ is a sequence of eigenmodes of the operator 
$\hat{P}_0(\hbar)$. In this case, the time evolution is trivial and one has $\mu_{\hbar}(t)=\mu_{\hbar}(0)$ for every $t$ in $\IR$. In other words, one only 
has to understand the distribution at time $t=0$. In the context of a chaotic Hamiltonian flow, a typical result concerning the properties of $\mu_{\hbar}$ is 
the Quantum Ergodicity property~\cite{Sh74, Ze87, CdV85, HeMaRo87} which we now recall.

Consider two sequences $(a_{\hbar})_{\hbar>0}$ and $(b_{\hbar})_{\hbar>0}$ such that $a_{\hbar},b_{\hbar}\rightarrow 1/2$ and $b_{\hbar}-a_{\hbar}\geq\alpha\hbar$ 
for some fixed positive constant $\alpha$. Then, take an orthonormal basis $(\psi_{\hbar}^j)_{j=1,\ldots, N(\hbar)}$ of 
$\mathbf{1}_{[a_{\hbar},b_{\hbar}]}(\hat{P}_0(\hbar))L^2(M)$ made of eigenmodes of $\hat{P}_0(\hbar)$. Suppose that the disintegration $L$ of the 
Liouville measure on $T_{1/2}^*M$ is ergodic for the geodesic flow, i.e.
$$\forall a\in\ml{C}^0(T_{1/2}^*M),\ \lim_{T\rightarrow+\infty}\frac{1}{T}\int_0^T a\circ G_0^t(\rho) dt=\int_{T_{1/2}^*M}adL\ \ \text{a.e.}.$$
The main examples of manifolds that enjoy this property are given by those which are negatively curved. Under this ergodicity assumption, one can prove~\cite{Sh74, Ze87, CdV85, HeMaRo87} the existence, for every $\hbar>0$, of a 
subset $J(\hbar)\subset\{1,\ldots N(\hbar)\}$ satisfying $\lim_{\hbar\rightarrow 0}\frac{\sharp J(\hbar)}{N(\hbar)}=1$, and
$$\forall a\in\ml{C}^{\infty}_c(T^*M), \lim_{\hbar\rightarrow 0, j\in J(\hbar)}\left\la \Oph(a)\psi_{\hbar}^j,\psi_{\hbar}^j\right\ra_{L^2(M)}=\int_{T_{1/2}^*M}adL.$$ 
In other words, the chaotic features of the classical dynamics imply that most of the stationary solutions of~\eqref{e:schrodinger} are equidistributed in the semiclassical limit. This result has many generalizations, with much progress made recently in terms of the understanding of the accumulation points along the exceptional subset $J(\hbar)^c$ -- we refer the reader to the following recent surveys~\cite{No13, Sa11, Ze10}.

A second natural class of initial data is given by a sequence of coherent states $(\psi_{\hbar}^{\rho_0})_{\hbar>0}$ which are microlocalized at a point $\rho_0$ 
in $T_{1/2}^*M$. In this case, one can use WKB methods to give a very precise description of the solution 
$u_{\hbar}^{\rho_0}(t)$ as long as $|t|\leq\kappa_0|\log\hbar|$ where $\kappa_0>0$ is some upper bound on the Lyapunov exponents of the geodesic flow~\cite{CoRo97}. In particular, under the ergodicity assumption of the Liouville measure and for any sequence $\tau_{\hbar}$, such that $\tau_{\hbar} \leq\kappa_0|\log\hbar|$, tending to $+\infty$, one can verify that, for a.e. choice of $\rho_0$ in $T_{1/2}^*M$ for any $0\leq c_1<c_2\leq 1$ and for any $a\in\ml{C}^{\infty}_c(T^*M)$, 
\begin{equation}\label{e:coherent}\lim_{\hbar\rightarrow 0}\int_{c_1}^{c_2}\left\la \Oph(a)u_{\hbar}^{\rho_0}(t\tau_{\hbar}),u_{\hbar}^{\rho_0}(t\tau_{\hbar})\right\ra_{L^2(M)}dt=(c_2-c_1)\int_{T_{1/2}^*M}adL.\end{equation} 
Another possibility is to consider a sequence of Lagrangian states $(\psi_{\hbar}^{\ml{L}})_{\hbar>0}$, which are microlocalized on a 
Lagrangian submanifold $\ml{L}$. In this case, it was proved that \emph{for a negatively curved manifold} and for any sequence 
$\tau_{\hbar}$ (again below $\kappa_0|\log\hbar|)$ tending to $+\infty$, property~\eqref{e:coherent} holds for a 
\emph{generic choice}\footnote{The precise assumption is that the Lagrangian submanifold must be transverse to the stable manifold on a set of large dimension.} 
of Lagrangian submanifolds~\cite{Sch05}. In this case, the equidistribution property is even stronger as it holds for any $0< t\leq 1$, i.e. one does not 
have to average in time. The equidisdribution properties for these two natural examples of initial data are in some sense non-stationary versions 
of Quantum Ergodicity. We also refer the reader to~\cite{BoDB00, BouDB05, AnNo07, Fa07, AnRi12} for related results on the description of the long time quantum dynamics
when the underlying classical system is chaotic.


A common feature of the above examples is that we can observe equidistribution in phase space of the solutions of~\eqref{e:schrodinger} under some chaoticity assumptions on the classical dynamics. We underline that, in all these examples, we are able to deduce  equidistribution because the initial data enjoy some additional properties. Yet, if we are given a general sequence of initial data $(\psi_{\hbar})_{\hbar>0}$ satisfying~\eqref{e:freq-assumption-unit}, we are a priori not able to say anything on the semiclassical behavior of the corresponding solutions, even for short logarithmic times.


\subsection{Perturbation of the Schr\"odinger equation}

In~\cite{Pe84}, motivated by the fact that the unitarity of~\eqref{e:schrodinger} does not allow to observe any sensitivity to the initial conditions in the quantum setting, Peres suggested a mechanism to observe ``irreversibility'' both in classical and quantum mechanics. His strategy is driven by the following principle: \emph{Instead of assuming that our preparations are marred by limited accuracy, we may assume that they are perfect but, on the other hand, the Hamiltonian $p_0$ is not exactly known, because we cannot perfectly insulate the physical system from its environment}~\cite{Pe84}. According to Peres, this point of view should allow to distinguish regular systems from chaotic ones both at the classical level and at the quantum level. In his article, he introduced quantum fidelity (which is now often called Loschmidt echo) as a way to measure the sensitivity of the quantum evolution to perturbations. In the last fifteen years, the study of this quantity has lead to a vast literature in the physics community and much progress has been made in the understanding of the quantum Loschmidt echo and its dependence to the properties of the underlying classical system. As this is the central object in the study of these questions in the physics literature, and as our article will be motivated by closely related questions, we briefly recall its definition translated into our context and some of its properties taken from the physics literature. We refer to the review articles~\cite{GPSZ06, JaPe09, GJPW12} for a detailed and precise account on these questions and for references to the existing literature.

In order to define the quantum Loschmidt echo, we fix a sequence of initial data $(\psi_{\hbar})_{\hbar>0}$ satisfying~\eqref{e:freq-assumption-unit}. We still denote by $u_{\hbar}(t)$ the solution at time $t$ of~\eqref{e:schrodinger} and we introduce $u_{\hbar}^{\eps}(t)$ the solution of the same equation where we replace $\hat{P}_0(\hbar)$ by 
$$\hat{P}_{\eps}(\hbar)=\hat{P}_0(\hbar)+\eps_0 V_0+\eps_1 V_1+\ldots+\eps_JV_J,$$ 
with $V_j$ smooth potentials on $M$. The quantum Loschmidt echo is then defined as the fidelity between these two evolved states, i.e.
$$F_{\hbar,\eps}(t):=\left|\left\la u_{\hbar}^{\eps}(t), u_{\hbar}(t)\right\ra_{L^2(M)}\right|^2.$$
According to the physics literature, this quantity is expected to decay for any quantum system, and the involved decay rate allows 
to distinguish chaotic sytems from regular ones. We will focus more precisely on the setting where the underlying dynamics is strongly chaotic, 
say the configuration space $M$ is a negatively curved manifold. In this particular context, we will now describe the properties that $F_{\hbar,\eps}(t)$ 
(or its average over a family of typical initial data, or perturbations) should satisfy due to the chaotic nature of the geodesic flow. We emphasize that this particular model has probably not been considered in the physics literature, and that we only translate in our setting the general properties which are expected for a typical strongly chaotic system.

In the case of a strong perturbation at the quantum level, i.e. which is large compared with the mean level spacing of the unperturbed 
system\footnote{In our setting, it is expected to be of order $\hbar^d$.}, the quantum Loschmidt echo should, after a short transition regime, 
typically exhibit an exponential decay of order $e^{-\Gamma(\eps) t}$, where $\Gamma(\eps)>0$ is some constant related to the physical system 
and to the strength of the perturbation~\cite{JalPas01, JaSiBe01}. An important feature of this exponential decay is that, for strong enough perturbations, the rate 
$\Gamma(\eps)$ should be \emph{independent} of the size of the perturbation and that it could be expressed only in terms of the Lyapunov exponents of the classical system.
After this ``exponential decay'', the quantum Loschmidt echo is expected to reach a long time saturation regime where it should take a value of order $\hbar^d$. For the matter of comparison, we mention that, in the case of integrable systems, the quantum Loschmidt echo should typically exhibit an algebraic decay. Again, we emphasize that this kind of decay rate is expected for perturbations which are strong at the quantum level even if they correspond to small classical perturbations. In the case of small quantum perturbations, the situation seems to be slightly more delicate, and it is typically expected that the quantum Loschmidt echo will follow a Gaussian exponential decay (with a rate depending on $\eps$) before it reaches a similar long time saturation regime. In both situations, the time scales appearing in the physics literature can be rather large, and sometimes much longer than the Ehrenfest time introduced in the beginning of the introduction.


\section{Statement of the main result}

In this article, we will not consider quantities which are directly related to the quantum Loschmidt echo. However, our general aim is to study questions which are also motivated by the general principle proposed in~\cite{Pe84} -- section $2$. More precisely, we will consider \emph{a fixed sequence of normalized initial data} $(\psi_{\hbar})_{\hbar>0}$ satisfying~\eqref{e:freq-assumption-unit} and we will try to understand what are the properties of the quantum evolution under \emph{typical perturbations of the Schr\"odinger equation}. Regarding the previous discussion on the quantum Loschmidt echo, we underline that we will look at strong quantum perturbations where $\|\eps\| \gg\sqrt{\hbar}$, and that our result will concern a regime of times for which the quantum Loschmidt echo is expected to exhibit an exponential decay with a 
rate independent of $\eps$.

These kinds of questions were recently considered in~\cite{EsTo12} by the first author and Toth in the context of magnetic perturbations of the Schr\"odinger equation -- see also~\cite{CanJaTo12} for generalizations of these results to metric perturbations. In these references, the authors studied the Schr\"odinger equation on general compact Riemannian manifolds and obtained average pointwise bounds on solutions for typical families of perturbations of the Schr\"odinger equation. In the following, we will also be interested in the behaviour of the solutions for typical perturbations of the Schr\"odinger equation but, instead of looking at their pointwise bounds, we will consider their distributions~\eqref{e:dist-phase-space} on phase space. For that study, we will restrict ourselves to scales of times where the classical-quantum correspondence is known to be valid, i.e. below the Ehrenfest time and our main focus will be on situations where {the underlying classical system is chaotic}. The two main issues in this article will be: 
\begin{itemize}
 \item how do small perturbations of the Hamiltonian $p_0$ alter the classical motion of a fixed particle (Theorem~\ref{t:dynamics})?
 \item how the phenomena we prove at the classical level can be transposed to the quantum level (Theorem~\ref{t:maintheo})?
\end{itemize}

For that purpose, we will consider one of the simplest models of a closed chaotic Hamiltonian system. Precisely, we let $M$ be a smooth compact boundaryless Riemmanian \emph{surface} of constant negative curvature $K\equiv -1$. In this case, it is known that the geodesic flow $G_0^t$ acting on the unit cotangent bundle $T_{1/2}^*M$ has the Anosov property~\cite{Ano67, KaHa, Rug07} -- see also section~\ref{s:geom-back} for a brief reminder. 


We will now introduce \emph{perturbations} of the Schr\"odinger operator $\hat{P}_0(\hbar)$ (and thus of the classical Hamiltonian $p_0$). For that purpose, we consider a family of potentials $V(\eps,.)$ on $M$ which is indexed by a multi-parameter $\eps=(\eps_0,\eps_1,\ldots,\eps_J)\in\mathbb{R}^{J+1}$. Precisely, we set
\begin{equation}\label{e:form-perturbation}
V(\eps,x):=\eps_0V_0(x)+\eps_1 V_1(x)+\eps_2 V_2(x)+\ldots+\eps_JV_J(x), 
\end{equation}
where $J\geq 0$ is some integer and $V_j$ are smooth \emph{real-valued} functions on $M$ for every $0\leq j\leq J$. Then, we will be interested in the following family of Schr\"odinger equations:
\begin{equation}\label{e:pert-schrodinger}\imath\hbar\frac{\partial u_{\hbar}^{\eps}}{\partial t}=\hat{P}_{\eps}(\hbar) u_{\hbar}^{\eps},\qquad u_{\hbar}^{\eps}\rceil_{t=0}=\psi_{\hbar}\in
L^{2}(M),
\end{equation}
where, for every $\eps=(\eps_0,\eps_1,\ldots,\eps_J)$, we define 
\begin{equation}
\hat{P}_{\eps}(\hbar):=\hat{P}_0(\hbar)+ V(\eps,x).
\end{equation}
Our objective will be to understand the behavior of the solutions $(u_{\hbar}^{\eps}(t))_{\hbar>0}$ of~\eqref{e:pert-schrodinger} for a fixed choice of initial data $(\psi_{\hbar})_{\hbar>0}$, and for a generic choice of $\eps$ in a box $(-\eps_{\hbar},\eps_{\hbar})^{J+1}$, with $\eps_{\hbar}\rightarrow 0^+$. 

\begin{rema}
We note that, in a first version of the present work~\cite{EsRi14}, we considered perturbation with only $1$ parameter of the form
$$\tilde{V}(\delta,x)=\delta V_0(x)+\delta^{1+\gamma} V_1(x)+\delta^{1+2\gamma} V_2(x)+\ldots+\delta^{1+\gamma J}V_J(x),$$
with $\gamma>0$ small enough and $\delta\in[0,\eps_{\hbar}]$. Modulo some technical issues, we obtained there similar results. Yet, we consider here the multi-parameter case which makes some aspects of the exposition simpler.
\end{rema}

In order to ensure that the perturbation is nontrivial, we need to impose some admissibility conditions on the potential. Fix $(x_0,\xi_0)$ in $T_{1/2}^*M$ and denote $(x(t),\xi(t))=G_0^t(x_0,\xi_0)$. We introduce the following quantity, for 
$W$ in $\ml{C}^{\infty}(M,\IR)$, 
\begin{equation}\label{e:admissibilityconstant}
\ml{L}_{x_0,\xi_0}(W):=\frac{1}{2}\int_{0}^{+\infty}g_{x(t)}^*\left(d_{x(t)}W,\xi^{\perp}(t)\right) 
e^{-t}dt,
\end{equation}
where $\xi^{\perp}$ is the unit vector directly orthogonal to $\xi$. The function $(x,\xi)\mapsto g^*_x(d_xW,\xi^{\perp})$ can be interpreted as the component 
in the unstable direction of the Hamiltonian vector field $X_W$ associated to $W$. Then, our admissibility condition on the family of potential $(V_j)_{j=0,\ldots, J}$ will read
\begin{equation}\label{e:admissibility}
\forall\ (x_0,\xi_0)\in T_{1/2}^*M,\ \exists 0\leq j\leq J,\ \ml{L}_{x_0,\xi_0}(V_j)\neq 0.
\end{equation}
\begin{rema}
This dynamical condition will appear naturally in our proof when we will differentiate the strong structural stability equation~\cite{dLLMM86}. We will provide examples satisfying these assumptions (with $J\geq 1$ a priori) in paragraph~\ref{ss:examples}. We underline that it is not apriori clear if one can find an example where $J=0$. The fact that we require that this dynamical condition holds for every $(x_0,\xi_0)$ comes from the fact that we want to consider \emph{any sequence} of initial data. For instance, if we had restricted ourselves to a sequence of coherent states microlocalized at a point $(x_0,\xi_0)$, then our proof should a priori work, modulo some extra work, by requiring only $\ml{L}_{x_0,\xi_0}(V_0)\neq 0$.
\end{rema}




As an application of our different results, we will obtain the following theorem on the equidistribution of a fixed sequence of initial data under a ``random perturbation'' of the Schr\"odinger equation.

\begin{theo}\label{t:appl} Let $M$ be a smooth compact boundaryless Riemannian surface of constant negative curvature $K\equiv -1$. Suppose 
$(V_j)_{j=0,\ldots J}$ satisfies~\eqref{e:admissibility}. Let $(\eps_{\hbar})_{\hbar>0}$ be a sequence satisfying $\eps_{\hbar}\longrightarrow 0$, and
$$\exists\ 0<\nu<\frac{1}{2}\ \text{such that}\ \eps_{\hbar}\geq \hbar^{\nu},$$
for $\hbar>0$ small enough. Set 
$$\tau_{\hbar}:=\log\left(\frac{1}{\eps_{\hbar}}\right).$$
Then, for every sequence of normalized states $(\psi_{\hbar})_{\hbar>0}$ satisfying~\eqref{e:freq-assumption-unit}, one can find $J(\hbar)\subset(-\eps_{\hbar},\eps_{\hbar})^{J+1}$ satisfying
$$\lim_{\hbar\rightarrow 0}\frac{\text{Leb}(J(\hbar))}{(2\eps_{\hbar})^{J+1}}=1,$$
and, for any $1\leq c_1\leq c_2\leq \min\{3/2,1/(2\nu)\}$ and for any $a\in\ml{C}^{\infty}_c(T^*M)$, 
$$\lim_{\hbar\rightarrow 0,\eps\in J(\hbar)}\int_{c_1}^{c_2}\left\la \Oph(a)u_{\hbar}^{\eps}(t\tau_{\hbar}),u_{\hbar}^{\eps}(t\tau_{\hbar})\right\ra_{L^2(M)}dt=(c_2-c_1)\int_{T_{1/2}^*M}adL,$$
where $u_{\hbar}^{\eps}(t')$ is the solution at time $t'$ of~\eqref{e:pert-schrodinger} with initial condition $\psi_{\hbar}$.


\end{theo}


Theorem~\ref{t:appl} shares some similarities with the results described in paragraph~\ref{ss:unperturbed} for the unperturbed Schr\"odinger equation in the sense that we also obtain some equidistribution property under the quantum evolution. The main ``improvement'' compared with the above results is that \emph{our theorem holds for any choice of initial data} $(\psi_{\hbar})_{\hbar>0}$ satisfying the frequency assumption~\eqref{e:freq-assumption-unit}, i.e. which are microlocalized on the unit cotangent bundle $T_{1/2}^*M$. \emph{However, one cannot ensure equidistribution for any Schr\"odinger equation but only for most perturbation of the equation}, where the ``random parameter'' we have to choose depends on the choice of initial data. For instance, if one chooses the initial data to be a family of coherent states microlocalized in $\rho_0\in T_{1/2}^*M$ belonging to a closed geodesic of $G_0^t$, the WKB methods from~\cite{CoRo97} apply. In particular, one can verify that, for every $\eps$ in a small enough neighborhood of $0$ (say $\eps\in[-\hbar^2,\hbar^2]^{J+1}$), the accumulation point we obtain is just the Lebesgue measure along the closed geodesic. This illustrates that we cannot hope to prove equidistribution for every choice of perturbation, but only for ``generic choices'' of $\eps\in(-\eps_{\hbar},\eps_{\hbar})^{J+1}$.

Concerning the allowed time scales $\tau_{\hbar}\rightarrow +\infty$, we underline that we are in a situation where the semiclassical approximation is valid, i.e. below the Ehrenfest time~\cite{BaGrPa99, BoRo02, Zw12}. In fact, we have to restrict ourselves to times $\tau_{\hbar}$ of order $|\log\eps_{\hbar}|$ and the relations between the involved parameters imply that $t\tau_{\hbar}\leq |\log\hbar|/2$, where $1/2$ is the best known constant for the Egorov Theorem in this geometric context~\cite{AnNo07, DyGu14}. Combining the remarks following corollary~\ref{c:maincoro} to the long time Egorov theorem~\ref{t:large-time-egorov3}, one can also verify that, for every $0\leq t<1$ and for every $\eps\in(-\eps_{\hbar},\eps_{\hbar})^{J+1}$, one has, as $\hbar\rightarrow 0$,
$$\left\la \Oph(a)u_{\hbar}^{\eps}(t|\log\eps_{\hbar}|),u_{\hbar}^{\eps}(t|\log\eps_{\hbar}|)\right\ra_{L^2(M)}
=\left\la \Oph(a)u_{\hbar}(t|\log\eps_{\hbar}|),u_{\hbar}(t|\log\eps_{\hbar}|)\right\ra_{L^2(M)}+o(1),$$
where $u_{\hbar}(t')$ is the solution at time $t'$ of~\eqref{e:schrodinger} with initial condition $\psi_{\hbar}$. Thus, $|\log\eps_{\hbar}|$ is really the critical time scale for which we can observe the effect of the perturbation. It would be of course very interesting to understand if similar results hold for longer times, e.g. $\tau_{\hbar}\gg|\log\eps_{\hbar}|$. The fact that our proof does not allow to go beyond $(1+1/2)|\log\eps_{\hbar}|$ will be related to the H\"older regularity of the conjugating homeomorphism in the strong structural stability theorem~\cite{Ano67, dLLMM86} which is a crucial element of our argument.

The above theorem is stated for \emph{surfaces of constant negative curvature}, and it is natural to ask whether the result can be extended to variable curvature or not. In fact, most of our proof will be valid for surfaces of variable curvature except for one step (see paragraph~\ref{ss:constant}) which requires to consider the case of constant curvature. In this particular case, we have a very precise description of the distribution of small pieces of unstable manifolds under the geodesic flow~\cite{Fu73, Mar77, Bu90}. These classical properties allow us to conclude in a simple manner. However, we believe that this restriction can be overcome (at least under some pinching condition on the curvature) by studying more precisely the analogous equidistribution properties in variable curvature by using for instance the dynamical tools from~\cite{Li04, BaLi12, FaTs13}. 
Another natural question would be to study the higher dimensional case. In this case, the situation would probably be slightly more delicate (even in constant curvature) as, among other difficulties, one would have to use perturbations by a family with more parameters in order to use equidistribution properties of the unstable manifold which is of dimension $d-1$.

Finally, it is also natural to compare this result with the case of integrable systems. In the particular case of Zoll manifolds, one cannot expect any equidistribution property 
as long as $\tau_{\hbar}\ll\eps_{\hbar}^{-1}$ as, for these scales of times, the accumulation point is always given by the average (over a period of the geodesic flow) 
of the semiclassical measure of the initial data~\cite{MacRiv15}.

\subsection{A brief outline of the proof}

One of the main issues of our article is to study the convergence of integrals of purely classical nature and our main statement on this question is contained in Theorem~\ref{t:dynamics}. Even if the statement of this theorem is quite natural, it seems to be a new result on the chaotic properties of the geodesic flow on negatively curved surfaces. This result on perturbations of the classical flow has an analogue for perturbations of the Schr\"odinger equation. The quantum version of this perturbation result is given by Theorem~\ref{t:maintheo} whose proof is based on the validity of the semiclassical tools below the Ehrenfest time.

More precisely, the two main steps of the proof can be divided as follows. First, we use more or less standard arguments of semiclassical analysis in order to reduce ourselves to a problem concerning only the classical properties of the Hamiltonian flow. This is done in section~\ref{s:average}. After making this reduction, we have to understand the asymptotic behaviour as $b\rightarrow 0$ and $T_0\rightarrow+\infty$ (simultaneously) of integrals of the following type:
$$I_{x_0,\xi_0}(b,T_0):=\frac{1}{(2b)^{J+1}}\int_{(-b,b)^{J+1}}a\circ G_{\eps}^{T_0}(x_0,\xi_0)d\eps,$$
where $(x_0,\xi_0)$ is some point in a neighborhood the unit cotangent bundle $T_{1/2}^*M$, $a$ is some smooth function and $G_{\eps}^t$ is 
the Hamiltonian flow associated to $p_{\eps}(x,\xi)=\frac{\|\xi\|^2}{2}+ V(\eps,x)$. Again, even if this kind of question seems quite natural, we were not 
able find any place in the literature where this question (or a related one) was considered. We will prove that, for $T_0\sim|\log b|$, this integral 
converges to $\int_{T_{1/2}^*M}adL$ where $L$ is the disintegration of the Liouville measure on $T_{1/2}^*M$. Our precise statement on this question is contained in Theorem~\ref{t:dynamics} -- see also corollary~\ref{c:maincoro}.
In order to prove this convergence result, we will proceed in two steps:
\begin{itemize}
\item we apply the strong structural stability theorem for Anosov flows to the flow $G_{\eps}^t$~\cite{Ano67, dLLMM86} and we study some geometric properties of the conjugating homeomorphism (section~\ref{s:struct-stab});
\item then, in section~\ref{s:mixing}, we conclude using unique ergodicity of the horocycle flow~\cite{Fu73, Mar77, Bu90}.
\end{itemize}
We also underline that an important aspect we have to deal with is the fact that we need to prove that \emph{the convergence of the integral $I_{x_0,\xi_0}(b,T_0)$ is uniform with respect to $(x_0,\xi_0)$}.

Finally, even if the problems under consideration are slightly different, we observe that some aspects of our article share some similarities with the proof on the distribution of ``typical'' Lagrangian states given in~\cite{Sch05}, or more generally with the classical proof of Quantum Ergodicity. Compared with these results, the main new ingredient in our approach is probably the use of strong structural stability, and the precise description of the geometric properties of the conjugating homeomorphism. This provides a way to transform our problem into a more standard problem of equidistribution of unstable manifolds.

Again, we emphasize that our proof is valid for general surfaces of variable curvature except for one step, namely the unique ergodicity of the horocycle flow 
for the Liouville measure in paragraph~\ref{ss:constant}. Yet, we will write our proof in the general context of surfaces of variable negative curvature as long as we can (as 
it does not really complexify our proof). 




\subsection{Organization of the article}

In section~\ref{s:geom-back}, we briefly recall some classical facts on Riemannian geometry and on Anosov geodesic flows. In section~\ref{s:average}, we give the proof of Theorem~\ref{t:appl} and show how it restricts to a problem of classical dynamics. Sections~\ref{s:struct-stab} and~\ref{s:mixing} are devoted to the proof of the dynamical questions raised in section~\ref{s:average}. Finally, in appendix~\ref{a:pdo}, we give a brief reminder on semiclassical analysis and in appendix~\ref{a:map-manifold}, we briefly describe the differentiable structure on manifolds of mappings.

\section{Geometric and dynamical preliminaries}\label{s:geom-back}

In this section, we briefly recall classical results on Riemannian geometry and Anosov geodesic flows that we will use at different stages of this article. Along the way, we 
fix some notations that we will use all along this article. We refer for instance the reader to~\cite{Bes78, Rug07, Web99} for more details.\\

\textbf{From this point and for the rest of the article, we will restrict ourselves to the case where $\dim M=2$} even if many of the properties described in this section can be 
adapted to higher dimensions.

\subsection{Almost complex, Riemannian and symplectic structures on $T^*M$} 

In this paragraph, we collect some standard facts on Riemannian geometry that can be found in more details in~\cite{Bes78} (Chap.~$1$) or~\cite{Web99} (Appendix $B$).

\subsubsection{Musical isomorphisms}

Recall that the Riemannian metric $g$ on $M$ induces two natural isomorphisms
$$\flat :T_xM\rightarrow T_x^*M,\ v\mapsto g_x(v,.),$$
and its inverse $\sharp: T_x^*M\rightarrow T_xM.$ This natural isomorphism induces a positive definite form on $T_x^*M$ for which these isomorphisms 
are in fact isometries. We denote by $g^*$ the corresponding metric. We will use the notation $g_x^*(\xi,\xi):=\|\xi\|_x^2$ for $(x,\xi)$ in $T^*M$ and we will often omit 
the subscript $x$ in order to alleviate notations.

\subsubsection{Horizontal and vertical subspaces} Let $\rho=(x,\xi)$ be an element in $T^*M$. Denote by $\pi: T^*M\rightarrow M$ the canonical projection $(x,\xi)\mapsto x$. We introduce the so-called \emph{vertical} subspace:
$$\ml{V}_{\rho}:=\text{Ker}(d_{\rho}\pi)\subset T_{\rho}T^*M.$$
The fiber $T_x^*M$ is a submanifold of $T^*M$ that contains the point $(x,\xi)$. The tangent space to this submanifold at point $(x,\xi)$ is the vertical subspace $\ml{V}_{\rho}$ and it can be canonically identified with $T_x^*M$. We will now define the connection map. For that purpose, we fix $Z$ in $T_{\rho}T^*M$ and $\rho(t)=(x(t),\xi(t))$ a smooth curve in $T^*M$ such that $\rho(0)=\rho$ and $\rho'(0)=Z$. The connection map $\mathcal{K}_{\rho}:T_{\rho}T^*M\mapsto T_x^*M$ is the following application:
$$\mathcal{K}_{\rho}(Z):=\frac{D}{dt}\rho'(0)=\nabla_{x'(0)}\xi(0),$$
where $\frac{D}{dt}\rho'(t)$ is the covariant derivative of $\rho(t)$ along the curve $x(t)$~\cite{GaHuLa} -- section II.B. One can verify that this quantity depends only on the initial velocity $Z$ of the curve and not on the curve, and that the map is linear. The \emph{horizontal} space is given by the kernel of this linear application, i.e.
$$\ml{H}_{\rho}:=\text{Ker}\ml{K}_{\rho}\subset T_{\rho}T^*M.$$
There exists a natural vector bundle isomoprhism between the pullback bundle $\pi^*(TM\oplus T^*M)\rightarrow T^*M$ and the canonical bundle $TT^*M\rightarrow T^*M$. The restriction of this isomorphism on the fibers above $\rho=(x,\xi)\in T^*M$ is given by
$$\theta(\rho):T_{\rho} T^*M\rightarrow (\rho, T_x M\oplus T_x^*M), Z\mapsto (y,\eta):=(d_{\rho}\pi Z,\ml{K}_{\rho}Z).$$ 
These coordinates $(y,\eta)$ will allow us to express easily the different structures on $T^*M$.

\begin{rema}\label{r:coordinate-chart} 

Let $(U,\phi)=(u^i)$ be local coordinates on $M$. One can lift in a canonical way these coordinates into coordinates $(u^i,v_j)$ on the cotangent space $T^*M$ 
by using $(\Phi,T^*U)$, where $\Phi(x,\xi)=(\phi(x),(d\phi(x)^T)^{-1}\xi)$. When we make a change of coordinates $(\tilde{u}^i,\tilde{v}_j)$, one can verify that $\tilde{v}$ is the image of $v$ under a linear transformation (depending only on the coordinate charts $(u^i)$ and $(\tilde{u}^i)$). A similar procedure can be performed to obtain coordinates $(u^i,v_j;y^k,w_l)$ on $TT^*M$. In this second step, the coordinate $(y^k)$ will satisfy the same property as $(v_j)$ but the coordinate $(w_l)$ will not. Precisely, if we consider another system of coordinates $(\tilde{u}^i,\tilde{v}_j;\tilde{y}^k,\tilde{w}_l)$ build on $(\tilde{U},\tilde{\phi})=(\tilde{u}^i)$, we will have that under the change of coordinate, $\tilde{y}$ can be expressed as the image of $y$ under a linear transformation (depending only on the coordinate charts $(u^i)$ and $(\tilde{u}^i)$) while $\tilde{w}$ will be the image of $w$ under a map involving $y$, $v$. However, if we set $\eta_l=w_l-\Gamma^q_{lp}(u)y^pv_q$ where $\Gamma^q_{lp}$ are the Christoffel symbols of the Levi-Civita connection~\cite{GaHuLa}, one can verify that $\tilde{\eta}$ is the image of $\eta$ under a linear map depending only on the coordinate charts $(u^i)$ and $(\tilde{u}^i)$. In local coordinates, the horizontal and vertical spaces can be written as
$$\ml{H}_{\rho}:=\{(u^i,v_j;y^k, \Gamma^q_{lp}(u)y^pv_q):y\in\IR^2\},$$
and
$$\ml{V}_{\rho}:=\{(u^i,v_j;0, w^l):w\in\IR^2\}.$$
Finally, the map $\theta(\rho)$ can be expressed as follows:
$$\theta(\rho):T_{\rho} T^*M\rightarrow (\rho, T_x M\oplus T_x^*M),\ (u^i,v_j;y^k,w_l)\mapsto (u^i,v_j;y^k,\eta_l=w_l-\Gamma^q_{lp}(u)y^pv_q).$$
 We refer the reader to appendix $B$ of~\cite{Web99} for more details on these coordinate changes.

\end{rema}

\subsubsection{Symplectic structure on $T^*M$} Recall that the canonical contact form on $T^*M$ is given by the following expression:
$$\forall\rho=(x,\xi)\in T^*M,\ \forall Z\in T_{\rho}T^*M,\ \alpha_{x,\xi}(Z)=\xi(d_{\rho}\pi(Z)).$$
The canonical symplectic form on $T^*M$ can then be defined as $\Omega=d\alpha$. Using our natural isomorphism, this symplectic form can be written as 
$$\forall Z_1\cong  (y_1,\eta_1)\in T_{\rho}T^*M,\ \forall Z_2\cong  (y_2,\eta_2)\in T_{\rho}T^*M,\ \Omega_{\rho}(Z_1,Z_2)=\eta_1(y_2)-\eta_2(y_1).$$

\subsubsection{Almost complex structure on $T^*M$} One can define the following map from $T_xM\oplus T_x^*M$ to itself:
$$\tilde{J}_x(y,\eta)=(\eta^{\sharp},-y^{\flat}).$$
This map induces an almost complex structure on $T_{\rho} T^*M$ through the isomorphism $\theta(\rho)$. We denote this almost complex structure by $J_{\rho}$.

\subsubsection{Riemannian metric on $T^*M$} 

The Sasaki metric $g^S$ on $T^*M$ is then defined as
$$g^S_{\rho}(Z_1,Z_2):=g_x^*(\ml{K}_{\rho}(Z_1),\ml{K}_{\rho}(Z_2))+g_x(d_{\rho}\pi(Z_1),d_{\rho}\pi(Z_2)).$$
This is a positive definite bilinear form on $T_{\rho}T^*M$. The important point is that \emph{this metric is compatible with the symplectic structure on $T^*M$ through the almost complex structure}. Precisely, one has, for every $(Z_1,Z_2)\in T_{\rho}T^*M\times T_{\rho} T^*M$,
$$g^S_{\rho}(Z_1,Z_2)=\Omega_{\rho}(Z_1,J_{\rho}Z_2).$$
In fact, using the natural isomorphism, one has
$$\Omega_{\rho}(Z_1,J_{\rho}Z_2)=\eta_1(\eta_2^{\sharp})+y_2^{\flat}(y_1)=g^*_x(\eta_1,\eta_2)+g_x(y_1,y_2).$$

\subsubsection{Expression of the Hamiltonian vector field}

Fix a smooth real-valued function $V$ on $M$. We recall the expression (in horizontal/vertical coordinates) of the Hamiltonian vector field associated to the function
$$p_{V}(x,\xi):=\frac{\|\xi\|^2_x}{2}+ V(x).$$
For that purpose, we use the local coordinates defined in remark~\ref{r:coordinate-chart} and a direct calculation gives us
$$d_{u,v}p_{V}(y,w)=\eta^lg^{lj}(u)v_j+\frac{\partial V}{\partial u^k}y^k,$$
where $\eta^l=w_l-\Gamma_{lp}^q(u)y^pv_q$. Thus, for any $Z\cong (y,\eta)\in T_{\rho}T^*M$, we have, using the natural coordinates induced by the map $\theta$,
$$d_{x,\xi}p_V.Z=g^*_x(\xi,\eta)+g_x((d_xV)^{\sharp},y).$$
We denote by $X_{p_{V}}\cong (y_0,\eta_0)$ the Hamiltonian vector field. One has $d_{\rho}p_V.Z=\Omega_{\rho}(Z,X_{{p}_V})$. Thus, $(y_0,\eta_0)=(\xi^{\sharp},-d_xV)$ and, using the natural isomorphism $\theta$, one has finally $X_{p_V}(x,\xi)=\theta(x,\xi)^{-1}(\xi^{\sharp},-d_xV)$.

\subsubsection{An orthogonal basis on $T_{\rho}T^*M$} Let $\rho=(x,\xi)$ be an element in $T^*M$. We want to define an orthonormal basis of $T_{\rho}T^*M$ adapted to the splitting 
described above. Thanks to the fact that the manifold $M$ is oriented with a Riemannian structure, one can define a notion of rotation by $\pi/2$ in every cotangent space 
$T_x^*M$, which is of dimension $2$, and thus there exists an unique $\xi^{\perp}$ such that $(\xi,\xi^{\perp})$ is a direct orthogonal basis with 
$\|\xi\|=\|\xi^{\perp}\|$. 
We use this to define an orthogonal basis of $\ml{V}_{\rho}$:
$$Y_0(\rho):=\left(\theta(\rho)\right)^{-1}(0,\xi),\ \text{and}\ U(\rho)=\left(\theta(\rho)\right)^{-1}(0,\xi^{\perp}).$$
Then, we can define an orthogonal basis of $\ml{H}_{\rho}$ as follows
$$X_0(\rho)=J_{\rho}Y_0(\rho),\ \text{and}\ W_{\rho}=J_{\rho}U(\rho).$$
The family $(X_0,W,-JX_0,-JW)$ forms a direct orthogonal basis of $T_{\rho}T^*M$. The tangent vector $X_0(\rho)$ is the geodesic vector field, $Y_0(\rho)=-JX_0(\rho)$ is the direction tangent to the submanifold $\IR\xi$ and all the vectors in this basis are of norm $\|\xi\|$. The expression of $X_{p_V}$ in this basis is given by
$$X_{p_V}(x,\xi)=X_0(x,\xi)+\frac{g_x^*(d_xV,\xi)}{\|\xi\|^2}JX_0(x,\xi)+\frac{g_x^*(d_xV,\xi^{\perp})}{\|\xi\|^2}JW(x,\xi).$$

\subsubsection{Restriction to energy layers}

Fix $E>0$ and denote by $T_E^* M$ the codimension $1$ submanifold $\{(x,\xi):\|\xi\|^2_x/2=E\}$. We introduce $\pi_E:T_E^*M\rightarrow M$ the canonical projection and we 
define the $1$-dimensional subspace
$$\forall\rho=(x,\xi)\in T_E^*M,\ V_{\rho}:=\text{Ker}(d_{\rho}\pi_E)\subset\ml{V}_{\rho}.$$
One can check that $V_{\rho}$ is in fact equal to $\IR J_{\rho}W_{\rho}$ and that
$$T_{\rho}T_E^*M=\ml{H}_{\rho}\oplus V_{\rho}.$$
We denote by $H_{\rho}$ the subspace of $\ml{H}_{\rho}$ generated by $W_{\rho}$, which is orthogonal to the geodesic vector field $X_0(\rho)$. We also define the tangent 
normal bundle:
$$N_{\rho}=H_{\rho}\oplus V_{\rho}.$$
\begin{rema}
 We underline that $g^S$ induces a Riemannian metric on $T_E^*M$ for every $E>0$.
\end{rema}

Finally, we introduce the map 
$$\Pi:(x,\xi)\in T^*M - M \mapsto(x,\xi/\|\xi\|)\in T_{1/2}^*M.$$
We can observe that this map induces a map between $T_{\rho}T^*M$ and $T_{\Pi(\rho)}T_{1/2}^*M$ which has the following action on the orthogonal basis described above:
$$X_0(\rho)\mapsto \sqrt{2p_0(\rho)}X_0\circ\Pi(\rho),\ W_{\rho}\mapsto \sqrt{2p_0(\rho)}W_{\Pi(\rho)},\ J_{\rho}W_{\rho}\mapsto J_{\Pi(\rho)}W_{\Pi(\rho)}, \text{and}
\ Y_0(\rho)\mapsto 0,$$
where $p_0(\rho):=\frac{\|\xi\|_x^2}{2}$ for $\rho=(x,\xi)$.

\subsection{Anosov assumption}
\label{ss:anosov}

Suppose now that $M$ is a \textbf{surface with negative sectional curvature}. In this case, it is known that the geodesic flow\footnote{In this case, $G_0^t$ is the 
Hamiltonian flow associated to $p_0$.} $G_0^t$ satisfies the Anosov property on $T_{1/2}^*M$~\cite{Ano67, KaHa, Rug07}. Precisely, it means that, for every 
$\rho=(x,\xi)$ in $T_{1/2}^*M$, there exists a $G_0^t$-invariant splitting
\begin{equation}\label{e:anosov}
T_{\rho}T_{1/2}^*M=\IR X_0(\rho)\oplus E^u(\rho)\oplus E^s(\rho),
\end{equation}
where $X_0(\rho)$ is the Hamiltonian vector field associated to $p_0(x,\xi)=\frac{\|\xi\|_x^2}{2}$, $E^u$ is the unstable direction and $E^s(\rho)$ is 
the stable direction~\cite{Ano67, KaHa, Rug07}. In our setting, the unstable and the stable spaces are $1$-dimensional subspaces of $T_{\rho}T_{1/2}^*M.$ These three subspaces are preserved under the geodesic flow and there exist  constants $C_0>0$ and $\gamma_0>0$ such that, for any $t\geq 0$, for any $v^s\in E^s(\rho)$ and any $v^u\in E^u(\rho)$,
$$\|d_{\rho}G_0^tv^s\|_{G^t_0(\rho)}\leq C_0e^{-\gamma_0 t}\|v^s\|_{\rho}\ \text{and}\ \|d_{\rho}G_0^{-t}v^u\|_{G^{-t}_0(\rho)}\leq C_0e^{-\gamma_0 t}\|v^u\|_{\rho},$$
where $\|.\|_w$ is the norm associated to the Sasaki metric restricted to the tangent space $T_{w} T_{1/2}^*M$.

Thanks to the fact that we are on a negatively curved surface, one can also prove that the stable (and unstable) foliation associated to this splitting~\cite{Rug07} are 
of class $\ml{C}^{2-\eta}$ for every $\eta>0$~\cite{PuHi75, Has94}.


We conclude this paragraph by a precise description of the unstable and stable spaces in terms of the horizontal and vertical subspaces we have defined in the previous 
paragraph. From this point, most of the statements are really specific to the fact that we are on a negatively curved \emph{surface}. We refer the reader to Chapter $3$ 
(Section $1$) of~\cite{Rug07} for a more detailed exposition.

Let $\rho=(x,\xi)$ be an element in $T_{1/2}^*M$. The unstable (resp. stable) space can be written as
$$E^u(\rho):=\left\{ (I-U^u(\rho) J_{\rho}) W:W\in H_{\rho}\right\},$$
resp.
$$E^s(\rho):=\left\{ (I-U^s(\rho) J_{\rho}) W:W\in H_{\rho}\right\},$$
where $U^u>0$ (resp. $U^s<0$) is called the unstable (resp. stable) Riccati solution~\cite{Rug07}. In the case where the curvature is constant and equal to $-1$, one has $U^u=-U^s=1$. Recall that we defined a ``canonical'' unit vector $W_{\rho}$ in $H_{\rho}$. We can introduce
$$X^u_{1/2}(\rho):=(I-U^u(\rho) J_{\rho}) W_{\rho}\ \text{and}\ 
X^s_{1/2}(\rho):=(I-U^s(\rho) J_{\rho}) W_{\rho}.$$
The family $(X_0(\rho),X^u_{1/2}(\rho),X^s_{1/2}(\rho))$ form a basis of vectors the tangent space $T_{\rho}T_{1/2}^*M$ which is associated to the Anosov 
decomposition~\eqref{e:anosov}. We underline that $X^u_{1/2}$ and $X_{1/2}^s$ are not unit vectors and that they are not a priori orthogonal. These vector fields 
allows us to define a certain parametrization of the unstable horocycle flow~\cite{Mar77} that we will denote by $(H_u^{\tau})_{\tau\in\IR}$ and that satisfies
$$\forall\rho\in T_{1/2}^*M,\ \frac{d}{d\tau}\left(H_u^{\tau}(\rho)\right)=X^u_{1/2}\circ H_u^{\tau}(\rho).$$
Still according to~\cite{Rug07}, one also has 
\begin{equation}\label{e:unstablejac}
\forall\tau\in\IR,\ \frac{\|d_{G_0^{-\tau}(\rho)}G_0^{\tau}.X^u_{1/2}\circ G_0^{-\tau}(\rho)\|_{\rho}}{\|X^u_{1/2}(\rho)\|_{\rho}}=e^{-\int_0^{-\tau}U^u\circ G_0^s(\rho)ds}.
\end{equation}
Using these definitions, we can write the generalization in variable curvature of the operator defined by~\eqref{e:admissibilityconstant}:
\begin{equation}\label{e:admissibility-operator}
\ml{L}_{x_0,\xi_0}(a):=\int_{0}^{+\infty}g_{x(t)}^*\left(d_{x(t)}a,\xi^{\perp}(t)\right) 
\frac{e^{-\int_0^{t}U^u\circ G_0^s(x_0,\xi_0)ds}}{U^u\circ G_0^t(x_0,\xi_0)-U^s\circ G_0^t(x_0,\xi_0)}dt, 
\end{equation}
where $a$ is a smooth function on $M$. Finally, we introduce the following quantity related to the expansion rates of the geodesic flow:
$$U_+:=\max\left\{\max_{\rho\in T_{1/2}^*M}\{U^u(\rho)\},\max_{\rho\in T_{1/2}^*M}\{-U^s(\rho)\}\right\},$$
and
$$U_-:=\min\left\{\min_{\rho\in T_{1/2}^*M}\{U^u(\rho)\},\min_{\rho\in T_{1/2}^*M}\{-U^s(\rho)\}\right\}.$$

\begin{rema}\label{r:bound-derivative-flow} One can observe that there exists $C_0>0$ such that, for every $t$ in $\IR$, one has
$$\sup_{\rho\in T_{1/2}^*M}\|d_{\rho}G_0^t\|\leq C_0 e^{|t|U_+},$$
 where $\|.\|$ denotes the operator norm induced by the Sasaki metric on $T_{1/2}^*M$.
\end{rema}

\subsection{Examples of admissible potential}\label{ss:examples}

In this paragraph, we will provide examples of potential satisfying~\eqref{e:admissibility}. This is in fact a consequence of the following lemma:
\begin{lemm}\label{l:examples}
There exists $R_0>0$ such that, for every smooth compact Riemannian surface $M$ of constant negative sectional curvature $K\equiv -1$ with injectivity radius $\geq R_0$, and, 
for every $(x_0,\xi_0)\in T_{1/2}^*M$, $\ml{L}_{x_0,\xi_0}(a)\neq 0$ for every $a\in\ml{C}^{\infty}(M,\IR)$ outside a fixed hyperplane depending on $(x_0,\xi_0)$. 
\end{lemm}

\begin{rema} By hyperplane, we just mean the kernel of a nonzero continuous linear form such as $\ml{L}_{x_0,\xi_0}$. The fact that there exist families $(V_j)_{j=0,\ldots J}$ satisfying~\eqref{e:admissibility} is then a combination of this lemma with the continuity of the map $(x_0,\xi_0)\mapsto\ml{L}_{x_0,\xi_0}(a)$ (see lemma~\ref{l:holder}) and the compactness of $T_{1/2}^*M$.
\end{rema}

\begin{proof} First, we observe that $a\mapsto \ml{L}_{x_0,\xi_0}(a)$ is a linear form and that $|\ml{L}_{x_0,\xi_0}(a)|\leq C\|a\|_{\ml{C}^1}$. In particular, 
$\ml{L}_{x_0,\xi_0}$ is a continuous linear form on the Fr\'echet space $\ml{C}^{\infty}(M)$ endowed with the semi-norms induced by the $\ml{C}^k$-topology. Thus, it 
remains to prove that there exists a function $a$ for which the linear form does not vanish.

Suppose for simplicity that the injectivity radius $R_{\text{inj}}$ of $M$ is $\geq 2$. Fix $(x_0,\xi_0)$ in $T_{1/2}^*M$. Our first observation is that, thanks to the above description 
of the unstable vector field and of $X_{p_V}$, one has for every $a$ in $\ml{C}^{\infty}(M,\IR)$,
\begin{equation}\label{e:alternative-adm-op}
 \ml{L}_{x_0,\xi_0}(a)=\frac{1}{2}\int_0^{+\infty}\frac{d}{d\tau}\left(a\circ \pi_{1/2}\circ H^{\tau}_u\circ G_0^t(x_0,\xi_0)\right)_{\tau=0}e^{-t}dt.
\end{equation}
We can define a natural chart in a neighborhood of $x_0$ as follows. Take $\delta>0$ small enough and define
$$f_{x_0,\xi_0}:(t,\tau)\in[-1/2,1]\times[-\delta,\delta]\mapsto \pi_{1/2}\circ H^{\tau}_u\circ G^t_0(x_0,\xi_0)\in M.$$
We will now construct a potential which is compactly supported in the image of $f_{x_0,\xi_0}$. For that purpose, we fix $\chi_1\geq 0$ in $\ml{C}^{\infty}_c((-1/2,1))$ 
such that $\int_{\IR}\chi_1(t)e^{-t}dt=2$, and $\chi_2\geq 0$ in $\ml{C}^{\infty}_c((-\delta,\delta))$ such that $\chi_2'(0)=1$. We observe that these functions do not depend on $M$. 
We set 
$$a\circ \pi_{1/2}\circ H^{\tau}\circ G^t_0(x_0,\xi_0)=\chi_1(t)\chi_2(\tau).$$
The function $a$ can also be considered as a function on the universal cover $\mathbb{H}^2$ of $M$. We can observe that $\|da\|$ is in 
fact bounded by a constant $C_{\mathbb{H}^2}(\chi_1,\chi_2)>0$ that depends only on $\chi_1$, $\chi_2$ and the canonical metric on $\mathbb{H}^2$, but not on $M$. 
Finally, if we use formula~\eqref{e:alternative-adm-op}, we can write that
$$\ml{L}_{x_0,\xi_0}(a)\geq 1-C_{\mathbb{H}^2}(\chi_1,\chi_2)e^{-R_{\text{inj}}/2},$$
which concludes the proof of the lemma.

\end{proof}

\section{Averaging over the perturbation}\label{s:average}

In this section, we will make a series of semiclassical reductions that allow to reduce the proof of Theorem~\ref{t:appl} to a question concerning the classical properties 
of the geodesic flow $G_0^t$. Our main objective will in fact be to understand the average over the box $(-\eps_{\hbar},\eps_{\hbar})^{J+1}$ of the semiclassical quantities in which we are interested. More precisely, we want to study the following quantity, as $\hbar\rightarrow 0$,
$$\mu_{\hbar}^{\eps_{\hbar}}(t)(a):=\frac{1}{(2\eps_{\hbar})^{J+1}}\int_{(-\eps_{\hbar},\eps_{\hbar})^{J+1}}\left\la u_{\hbar}^{\eps}(t\tau_{\hbar}),
\Oph(a)u_{\hbar}^{\eps}(t\tau_{\hbar})\right\ra d\eps,$$
where $\tau_{\hbar}:=|\log(\eps_{\hbar})|$, $u_{\hbar}^{\eps}(t')$ is the solution at time $t'$ of~\eqref{e:pert-schrodinger} with initial condition $\psi_{\hbar}$, 
and $(\psi_{\hbar})_{\hbar>0}$ is a sequence of normalized states in $L^2(M)$ satisfying
\begin{equation}\label{e:freq-assumption-close-to-unit}
\mathbf{1}_{[1/2-\delta_0,1/2+\delta_0]}\left(\hat{P}_0(\hbar)\right)\psi_{\hbar}=\psi_{\hbar}+o(1),
\end{equation}
for some fixed $\delta_0\in(0,1/2)$. For simplicity of exposition, we will restrict ourselves to observables $a$ which are $0$-homogeneous in a small neighborhood of the 
unit cotangent bundle $T_{1/2}^*M$. Precisely, we fix $0\leq \tilde{\chi}\leq 1$ a smooth cutoff function on $\IR$ which is equal to $1$ in a small neighborhood of $1/2$ and which is 
equal to $0$ outside the interval $[1/4,3/4]$. Then, for every $a$ in $\ml{C}^{\infty}(T_{1/2}^*M)$, we define the following function which is compactly supported on 
$T^*M$:
$$\forall (x,\xi)\in T^*M,\ a_{\text{hom}}(x,\xi):=\tilde{\chi}\circ p_0(x,\xi)\times a\left(x,\frac{\xi}{\|\xi\|}\right).$$
With these conventions, we can state the following theorem which is in fact the main result of this section:

\begin{theo}\label{t:maintheo}  Let $M$ be a smooth compact boundaryless Riemannian surface of constant negative curvature $K\equiv -1$. Suppose 
$(V_j)_{j=0,\ldots J}$ satisfies~\eqref{e:admissibility}. Let $(\eps_{\hbar})_{\hbar>0}$ be a sequence satisfying, for $\hbar>0$ small enough,
$$\eps_{\hbar}\longrightarrow 0,\ \text{and}\ \exists\ 0<\nu<\frac{1}{2}\ \text{such that}\ \eps_{\hbar}\geq \hbar^{\nu}.$$

Then, for every $1<c_1<c_2<\min\{3/2,1/(2\nu)\}$, there exists $\delta_1>0$ such that, for any sequence of normalized states $(\psi_{\hbar})_{\hbar>0}$ satisfying~
\eqref{e:freq-assumption-close-to-unit} with $\delta_0=\delta_1$, and for any observable $a$ in $\ml{C}^{\infty}(T_{1/2}^*M)$, one has, uniformly for $t\in [c_1, c_2]$,
$$\lim_{\hbar\rightarrow 0}\frac{1}{(2\eps_{\hbar})^{J+1}}\int_{(-\eps_{\hbar},\eps_{\hbar})^{J+1}}\left\la u_{\hbar}^{\eps}(t|\log(\eps_{\hbar})|),
\Oph(a_{\operatorname{hom}})u_{\hbar}^{\eps}(t|\log(\eps_{\hbar})|)\right\ra d\eps=\int_{T_{1/2}^*M}adL,$$
where $L$ is the normalized Liouville measure on the energy layer $T_{1/2}^*M$ and $u_{\hbar}^{\eps}(t')$ is the solution at time $t'$ of~\eqref{e:pert-schrodinger} with initial condition $\psi_{\hbar}$.

\end{theo}

We emphasize that the range of times in this theorem correspond to times for which the semiclassical tools remain valid, i.e. below the Ehrenfest time. The proof of this theorem will be reduced to a question concerning the ergodic properties of the geodesic flow -- see paragraph~\ref{ss:reduction}. We will treat this question of purely classical nature in sections~\ref{s:struct-stab} and~\ref{s:mixing}.

It is interesting to mention that the formulation of this quantum result is very close to the formulation of its classical analogue, namely corollary~\ref{c:maincoro}. In both settings, we consider any fixed sequence of initial data localized near the unit cotangent bundle and we prove that the evolution enjoys some equidistribution properties when we average over a family of Hamiltonian close to $p_0$.  This theorem will allow us to prove Theorem~\ref{t:appl} -- see paragraph~\ref{ss:implication} below. For that purpose, one can remark that the previous theorem can be interpreted as a sort of local Weyl law which is usually a key argument in quantum ergodic statements.

\subsection{From Theorem~\ref{t:maintheo} to Theorem~\ref{t:appl}}\label{ss:implication}

The proof of our main result Theorem \ref{t:appl} follows a standard quantum ergodicity argument, which can be found in \cite{Sh74, Ze87, CdV85, HeMaRo87} (see also Chapter 
15 of \cite{Zw12}). The main differences with these references are that we consider a non-stationary problem, as in~\cite{AnRi12} for instance, and that 
we are averaging over a family of operators instead of a family of initial data. As in the above references, the proof can be divided in two main steps: 1) we estimate the variance, and 2) we extract a density $1$ subsequence.

We suppose that $\eps_{\hbar}$ satisfies the assumptions of theorem~\ref{t:appl}, and we fix a sequence of normalized states $(\psi_{\hbar})_{\hbar>0}$ in $L^2(M)$ 
which satisfies the frequency assumption~\eqref{e:freq-assumption-unit}. We are interested in the asymptotic behaviour, as $\hbar\rightarrow 0^+$ and 
$\eps\in(-\eps_{\hbar},\eps_{\hbar})^{J+1}$, of
$$\mu_{\hbar}^{\eps}\left(\mathbf{1}_{[c_1,c_2]}\otimes a\right):=\int_{c_1}^{c_2}\left\la \Oph(a)u_{\hbar}^{\eps}(t\tau_{\hbar}),u_{\hbar}^{\eps}(t\tau_{\hbar})\right\ra_{L^2(M)}dt,$$
when $1<c_1\leq c_2<\min\{3/2,1/(2\nu)\}$, $a$ in $\ml{C}^{\infty}_c(T^*M)$, and $\tau_{\hbar}:=|\log(\eps_{\hbar})|$.

\begin{rema}\label{r:adm-obs} Let $a$ be an element in $\ml{C}^{\infty}_c(T^*M)$. As above, we introduce
$$a_{\text{hom}}:=\tilde{\chi}\circ p_0(x,\xi)\times a\left(x,\frac{\xi}{\|\xi\|}\right).$$
Thanks to the frequency assumption~\eqref{e:freq-assumption-unit}, one can verify that
$$\limsup_{\hbar\rightarrow 0,\eps\in(-\eps_{\hbar},\eps_{\hbar})^{J+1}}\left|\mu_{\hbar}^{\eps}\left(\mathbf{1}_{[c_1,c_2]}\otimes a\right)-
\mu_{\hbar}^{\eps}\left(\mathbf{1}_{[c_1,c_2]}\otimes a_{\text{hom}}\right)\right|=0.$$
In particular, we can reduce our analysis to observables of the type $a_{\text{hom}}$ which are build from a smooth function $a$ on $T_{1/2}^*M$. 
\end{rema}

\subsubsection{Estimating the variance}

Let $a$ be an element in $\ml{C}^{\infty}(T_{1/2}^*M)$ and let $1<c_1\leq c_2<\min\{3/2,1/(2\nu)\}$. 

Our first step will be to establish the variance estimate in our setting by interpreting Theorem \ref{t:maintheo} as a sort of local Weyl law. More specifically, we will prove that 
\begin{align} \label{e:variance}
\frac{1}{(2\eps_{\hbar})^{J+1}}\int_{(-\eps_{\hbar},\eps_{\hbar})^{J+1}}\left| \mu_{\hbar}^{\eps}\left(\mathbf{1}_{[c_1,c_2]}\otimes a_{\text{hom}}\right) - (c_2-c_1)\int_{T_{1/2}^*M}adL \right|^2 d\eps \rightarrow 0
\end{align}
as $\hbar \rightarrow 0$. Thanks to the frequency assumption~\eqref{e:freq-assumption-unit}, we can assume, without loss of generality, that 
$\int_{T_{1/2}^*M}adL = 0$.  An important estimate that we will use is the following: for a fixed time $T_0$, one has 
$$d_{T^*M}(G_{\eps}^s (\rho), G_{0}^s(\rho))= \mathcal{O}(\|\eps\|),$$ 
uniformly for $\rho$ satisfying $p_0(\rho)\in[1/4,3/4]$, for $0\leq s\leq T_0$, and for $\eps\in(-1,1)^{J+1}$. Consequently, we then have, uniformly for $0\leq s\leq T_0$ and for $\eps\in(-\eps_{\hbar},\eps_{\hbar})^{J+1}$,
\begin{align*}
\mu_{\hbar}^{\eps}\left(\mathbf{1}_{[c_1,c_2]}\otimes a_{\text{hom}}\right) & = \int_{c_1}^{c_2} \left\la \Oph(a_{\text{hom}})u_{\hbar}^{\eps}(t\tau_{\hbar}),u_{\hbar}^{\eps}(t\tau_{\hbar})\right\ra_{L^2(M)}dt \\
& = \frac{1}{\tau_{\hbar}}\int_{c_1\tau_{\hbar}}^{c_2\tau_{\hbar}} \left\la  \Oph(a_{\text{hom}} \circ G_{\eps}^s)e^{\frac{is}{\hbar}\hat{P}_{\eps}(\hbar)} u_{\hbar}^{\eps}(t),e^{\frac{is}{\hbar}\hat{P}_{\eps}(\hbar)}u_{\hbar}^{\eps}(t)\right\ra_{L^2(M)}dt+  \mathcal{O}(\hbar) \\
& =  \frac{1}{\tau_{\hbar}}\int_{c_1\tau_{\hbar}}^{c_2\tau_{\hbar}} \left\la \Oph(a_{\text{hom}} \circ G_{0}^s)e^{\frac{is}{\hbar}\hat{P}_{\eps}(\hbar)} u_{\hbar}^{\eps}(t),e^{\frac{is}{\hbar}\hat{P}_{\eps}(\hbar)}u_{\hbar}^{\eps}(t)\right\ra_{L^2(M)}dt  + \mathcal{O}(\eps_{\hbar}) +  \mathcal{O}(\hbar) ,
\end{align*}
where we have first applied the Egorov theorem \ref{t:large-time-egorov3} for the fixed time parameter $0\leq s\leq T_0$ (with $\overline{\nu}=0$) and then used that 
$|a_{\text{hom}} \circ G_{0}^s(\rho) - a_{\text{hom}} \circ G_{\eps}^s(\rho) | = \mathcal{O}(\eps_{\hbar})$.  The property that $e^{\frac{is}{\hbar}\hat{P}_{\eps}(\hbar)}u_{\hbar}^{\eps}(t) = u_{\hbar}^{\eps}(t-s)$, 
along with a linear change of variables in $t$ give us the quantity
$$
\mu_{\hbar}^{\eps}\left(\mathbf{1}_{[c_1,c_2]}\otimes a_{\text{hom}}\right) = \frac{1}{\tau_{\hbar}}\int_{c_1\tau_{\hbar}-s}^{c_2\tau_{\hbar}-s}  \left\la \Oph(a_{\text{hom}} \circ G_{0}^s) u_{\hbar}^{\eps}(t),u_{\hbar}^{\eps}(t)\right\ra_{L^2(M)} dt+ \mathcal{O}(\eps_{\hbar}) +  \mathcal{O}(\hbar) .
$$
For our purposes, it is necessary to estimate away certain portions of the average in $t$. We rewrite the above quantity as 
\begin{align*}
\mu_{\hbar}^{\eps}\left(\mathbf{1}_{[c_1,c_2]}\otimes a_{\text{hom}}\right)  & = \frac{1}{\tau_{\hbar}}\int_{c_1\tau_{\hbar}}^{c_2\tau_{\hbar}}  \left\la \Oph(a_{\text{hom}} \circ G_{0}^s) u_{\hbar}^{\eps}(t),u_{\hbar}^{\eps}(t)\right\ra_{L^2(M)}  dt \\  
 & - \frac{1}{\tau_{\hbar}}\int_{c_2\tau_{\hbar}-s}^{c_2\tau_{\hbar}}  \left\la \Oph(a_{\text{hom}} \circ G_{0}^s )u_{\hbar}^{\eps}(t),u_{\hbar}^{\eps}(t)\right\ra_{L^2(M)}  dt \\
 & + \frac{1}{\tau_{\hbar}}\int_{c_1\tau_{\hbar}-s}^{c_1\tau_{\hbar} }  \left\la \Oph(a_{\text{hom}} \circ G_{0}^s) u_{\hbar}^{\eps}(t),u_{\hbar}^{\eps}(t)\right\ra_{L^2(M)}  dt \\
 & + \mathcal{O}(\eps_{\hbar}) +  \mathcal{O}(\hbar).
\end{align*}
It is straightforward to see that the second and third integrals are $\mathcal{O}(\frac{s}{\tau_h})$ after using the unitary of $e^{-\frac{it}{\hbar}\hat{P}_{\eps}(\hbar)}$ on $L^2$ along with the 
$L^2$ boundedness of zeroth order pseudodifferential operators by the Calderon-Vaillancourt Theorem. Furthermore, as $s$ is an artificial parameter, we can perform 
the average $ \la a_{\text{hom}} \ra_{T_0}  (\rho)= 1/{T_0} \int_{0}^{T_0}  a_{\text{hom}} \circ G_0^s(\rho) ds$ (for $\rho \in T^*M$) and $T_0>0$ the allowable time scale 
described in Theorem \ref{t:large-time-egorov3} (with $\overline{\nu}=0$).  Thus, performing the reverse change of variables $t = t' \tau_h$ in order to place us back into the framework of Theorem \ref{t:maintheo}, we finally arrive at 
$$
\mu_{\hbar}^{\eps}\left(\mathbf{1}_{[c_1,c_2]}\otimes a_{\text{hom}}\right)  = \mu_{\hbar}^{\eps}\left(\mathbf{1}_{[c_1,c_2]}\otimes \la a_{\text{hom}}\ra_{T_0}\right) 
+\mathcal{O}(\eps_{\hbar}) +  \mathcal{O}(\hbar) +\ml{O}(\tau_{\hbar}^{-1}),
$$
where the constants in the remainder depends on $a_{\text{hom}}$, $c_1$, $c_2$ and $T_0$. Recalling that we want to show the convergence of 
$\frac{1}{\eps_{\hbar}}\int_{0}^{\eps_{\hbar}}| \mu_{\hbar}^{\eps}\left(\mathbf{1}_{[c_1,c_2]}\otimes a_{\text{hom}}\right) |^2 d\eps,$ an application Jensen's inequality for convex functions on the $t$-integral shows
\begin{align*}
\frac{1}{(2\eps_{\hbar})^{J+1}}\int_{(-\eps_{\hbar},\eps_{\hbar})^{J+1}} & | \mu_{\hbar}^{\eps}\left(\mathbf{1}_{[c_1,c_2]}\otimes a_{\text{hom}}\right) |^2 d\eps  \\
& \leq \frac{(c_2-c_1)}{(2\eps_{\hbar})^{J+1}} \int_{(-\eps_{\hbar},\eps_{\hbar})^{J+1} }  \bigg[  \int_{c_1 }^{c_2 } 
\bigg| \left\la \Oph(\la a_{\text{hom}} \ra_{T_0}) u_{\hbar}^{\eps}( t\tau_{\hbar}),u_{\hbar}^{\eps}(t\tau_{\hbar})\right\ra_{L^2(M)}  \bigg|^2 \, dt \bigg] \, d\eps \\
& +\mathcal{O}(\eps_{\hbar}) +  \mathcal{O}(\hbar) +\ml{O}(\tau_{\hbar}^{-1}) \\
& =   \int_{c_1 }^{c_2 }\bigg[   \frac{(c_2-c_1)}{(2\eps_{\hbar})^{J+1}} \int_{(-\eps_{\hbar},\eps_{\hbar})^{J+1} } 
\bigg| \left\la \Oph(\la a_{\text{hom}} \ra_{T_0}) u_{\hbar}^{\eps}( t\tau_{\hbar}),u_{\hbar}^{\eps}(t\tau_{\hbar})\right\ra_{L^2(M)}  \bigg|^2 \, d\eps \bigg] \, dt \\
& +\mathcal{O}(\eps_{\hbar}) +  \mathcal{O}(\hbar) +\ml{O}(\tau_{\hbar}^{-1}).
\end{align*}
Furthermore, 
\begin{align*}
|\left\la \Oph(\la a_{\text{hom}} \ra_{T_0}) u_{\hbar}^{\eps}( t\tau_{\hbar}),u_{\hbar}^{\eps}(t\tau_{\hbar})\right\ra_{L^2(M)}|^2 \\ 
& \leq  \la \Oph(\la a_{\text{hom}} \ra_{T_0})^* \Oph(\la a_{\text{hom}} \ra_{T_0}) u_{\hbar}^{\eps}(t\tau_{\hbar}),u_{\hbar}^{\eps}(t\tau_{\hbar})\ra_{L^2(M)}.
\end{align*}
Then, using composition of pseudodifferential operators, one finds
\begin{align*}
\frac{1}{(2\eps_{\hbar})^{J+1}} \int_{(-\eps_{\hbar},\eps_{\hbar})^{J+1} } & | \mu_{\hbar}^{\eps}\left(\mathbf{1}_{[c_1,c_2]}\otimes a_{\text{hom}}\right) |^2 d\eps  \\
& \leq  \int_{c_1 }^{c_2} \left[\frac{(c_2-c_1)}{(2\eps_{\hbar})^{J+1}} \int_{(-\eps_{\hbar},\eps_{\hbar})^{J+1} } \la \Oph(|\la a_{\text{hom}} \ra_{T_0}|^2) u_{\hbar}^{\eps}(t\tau_{\hbar}),u_{\hbar}^{\eps}(t\tau_{\hbar})\ra_{L^2(M)}\, d\eps\right] dt\\
& +\mathcal{O}(\eps_{\hbar}) +  \mathcal{O}(\hbar) +\ml{O}(\tau_{\hbar}^{-1}).
\end{align*}
Observe now that $\la a_{\text{hom}} \ra_{T_0}=(\la a \ra_{T_0})_{\text{hom}}$. In particular, one finds that $|\la a_{\text{hom}} \ra_{T_0}|^2\leq 
(|\la a \ra_{T_0}|^2)_{\text{hom}}$. Then, we apply the positivity properties of paragraph~\ref{ss:antiwick} to write that
\begin{align*}
\frac{1}{(2\eps_{\hbar})^{J+1}} \int_{(-\eps_{\hbar},\eps_{\hbar})^{J+1} } & | \mu_{\hbar}^{\eps}\left(\mathbf{1}_{[c_1,c_2]}\otimes a_{\text{hom}}\right) |^2 d\eps  \\
& \leq  \int_{c_1 }^{c_2} \left[\frac{(c_2-c_1)}{(2\eps_{\hbar})^{J+1}} \int_{(-\eps_{\hbar},\eps_{\hbar})^{J+1} } \la \Oph((|\la a \ra_{T_0}|^2)_{\text{hom}}) u_{\hbar}^{\eps}(t\tau_{\hbar}),u_{\hbar}^{\eps}(t\tau_{\hbar})\ra_{L^2(M)}\, d\eps\right] dt\\
& +\mathcal{O}(\eps_{\hbar}) +  \mathcal{O}(\hbar) +\ml{O}(\tau_{\hbar}^{-1}).
\end{align*}
\begin{rema}
 We underline that, up to this point, our proof is valid for any sequences $\eps_{\hbar}\rightarrow 0$ and $\tau_{\hbar}\rightarrow+\infty$. In order to conclude, we will 
need to specify some restrictions on these sequences, namely the restrictions of Theorem~\ref{t:maintheo}.
\end{rema}

As the initial data $(\psi_{\hbar})_{\hbar>0}$ satisfy the frequency assumptions~\eqref{e:freq-assumption-unit}, one can finally apply Theorem~\ref{t:maintheo}. Keeping in mind that the convergence in this theorem is uniform for 
$t \in [c_1,c_2]$, we deduce that
$$
\limsup_{\hbar \rightarrow 0} \frac{1}{(2\eps_{\hbar})^{J+1}} \int_{(-\eps_{\hbar},\eps_{\hbar})^{J+1} } | \mu_{\hbar}^{\eps}\left(\mathbf{1}_{[c_1,c_2]}\otimes a_{\text{hom}}\right) |^2 d\eps 
 \leq (c_2-c_1)^2\int_{T^*_{1/2}M} \big| \la a \ra_{T_0} \big|^2 \, dL.
$$
This is valid for any $T_0>0$. As the Liouville measure is ergodic for $G_0^t$ in our geometric setting, we can conclude by letting $T_0\rightarrow +\infty$ that
$$
\limsup_{\hbar \rightarrow 0} \frac{1}{(2\eps_{\hbar})^{J+1}} \int_{(-\eps_{\hbar},\eps_{\hbar})^{J+1} } | \mu_{\hbar}^{\eps}\left(\mathbf{1}_{[c_1,c_2]}\otimes a_{\text{hom}}\right) |^2 d\eps=0,$$
which concludes our proof of the variance estimate.

\subsubsection{Extraction of a density $1$ subsequence}

Once we have the variance estimate~\eqref{e:variance}, the proof of Theorem~\ref{t:appl} follows directly after the standards arguments for extracting density 1 sequences as presented in, for example, \cite{Ze87, CdV85, HeMaRo87, Zw12}. As this argument is quite standard, we will just give a brief outline of the proof, refering to~\cite{Zw12} (Theorem $15.5$) for more details.


First, fix an observable $a \in C^{\infty}(T_{1/2}^*M)$ and $1<c_1<c_2<\min\{3/2, 1/(2\nu)\}.$ An application of Chebyshev's inequality combined with the variance estimate~\eqref{e:variance} allows to construct, for every $0<\hbar\leq 1$, $\Lambda(\hbar)\subset(-\eps_{\hbar},\eps_{\hbar})^{J+1}$ such that $\frac{|\Lambda(\hbar)|}{(2\eps_{\hbar})^{J+1}} \rightarrow 1$ as $\hbar \rightarrow 0$, and
$$\lim_{\hbar\rightarrow 0,\eps\in\Lambda(\hbar)}\mu_{\hbar}^{\eps}\left(\mathbf{1}_{[c_1,c_2]}\otimes a_{\text{hom}}\right) =(c_2-c_1)\int_{T_{1/2}^*M}adL.$$


Consider now $\{ a^k \}_{k=1}^{+\infty} \subset C^{\infty}(T_{1/2}^*M)$ which is dense in $\ml{C}^{0}(T_{1/2}^*M)$ for the $\ml{C}^0$ topology. Denote also by $\{c_l\}_{l=1}^{+\infty}$ the set of rational numbers in the open interval $(1,\min\{3/2, 1/(2\nu)\})$. Arguing as in~\cite{Zw12}, one can construct $\Lambda_{\infty}(\hbar)$ of density $1$ such that, for every $(k,l_1,l_2)\in\IN^3$ satisfying $c_{l_1}\leq c_{l_2}$, one has
$$\lim_{\hbar\rightarrow 0,\eps\in\Lambda_{\infty}(\hbar)}\mu_{\hbar}^{\eps}\left(\mathbf{1}_{[c_{l_1},c_{l_2}]}\otimes a^k_{\text{hom}}\right) =(c_{l_2}-c_{l_1})\int_{T_{1/2}^*M}a^kdL.$$


We set $J(\hbar) = \Lambda_{\infty}(\hbar)$, and we will now check that Theorem~\ref{t:appl} holds with this choice. We fix $a$ in $\ml{C}^{\infty}_c(T^*M)$ and $1\leq c_1<c_2\leq \min\{3/2, 1/(2\nu)\}$. Thanks to remark~\ref{r:adm-obs}, the accumulation points depend only on the value of $a$ on the unit cotangent bundle. Thus, we can without loss of generality replace $a$ by $a_{\text{hom}}$ (which coincides with $a$ on $T_{1/2}^*M$).

Let $0<\eta<1$. By density of the family $\{ a^k \}_{k=1}^{+\infty}$ in the $\ml{C}^0$ topology, there exists $k$ and $\hbar_0$ such that for all $\hbar < \hbar_0$
\begin{align*}
\int_{T^*_{1/2}M} |a^k - a_{\text{hom}}| \, dL < \eta \text{ and } \|\Oph(a_{\text{hom}}^k) - \Oph(a_{\text{hom}}) \|_{L^2 \rightarrow L^2} < \eta.
\end{align*}
One can also find $l_1$ and $l_2$ such that $|c_1-c_{l_1}|\leq\eta$ and $|c_2-c_{l_2}|\leq\eta$. Approximating by these different test functions, one finds that,
$$\limsup_{\hbar\rightarrow 0,\eps\in\Lambda_{\infty}(\hbar)}\left|\mu_{\hbar}^{\eps}\left(\mathbf{1}_{[c_{1},c_{2}]}\otimes a\right) -(c_{2}-c_{1})\int_{T_{1/2}^*M}adL \right|\leq C\eta,$$
where $C>0$ is some positive constant that does not depend on $\eta$. As this property holds for any $0<\eta<1$, the proof of Theorem~\ref{t:appl} is completed.


\subsection{Reduction of the proof of Theorem~\ref{t:maintheo} to a dynamical systems question}
\label{ss:reduction}

The purpose of this paragraph is to use the semiclassical approximation in order to obtain lower and upper bounds on 
$$\mu_{\hbar}^{\eps_{\hbar}}(t)(a)=\frac{1}{(2\eps_{\hbar})^{J+1}} \int_{(-\eps_{\hbar},\eps_{\hbar})^{J+1} }\left\la u_{\hbar}^{\eps}(t\tau_{\hbar}),\Oph(a)
u_{\hbar}^{\eps}(t\tau_{\hbar})\right\ra d\eps,$$
where $u_{\hbar}^{\eps}(t')$ is the solution at time $t'$ of~\eqref{e:pert-schrodinger} with initial data given by a sequence of \emph{normalized states} $(\psi_{\hbar})_{\hbar}$ satisfying the following frequency assumption:
\begin{equation}\label{e:energy-loc-states}\mathbf{1}_{[E_1,E_2]}(\hat{P}_0(\hbar))\psi_{\hbar}=\psi_{\hbar}+o(1),
\end{equation}
with $0<E_1\leq E_2<+\infty$. Then, the proof of theorem~\ref{t:maintheo} will be reduced to the computation of these bounds as $\hbar$ goes to $0$.

Precisely, we introduce the following quantities, for $\delta>0$ (small enough),
$$A_-(\hbar,t,\delta):=\inf_{E\in[E_1-\delta/2,E_2+\delta/2]}\inf_{\rho\in T_E^*M}\left\{\frac{1}{(2\eps_{\hbar})^{J+1}} \int_{(-\eps_{\hbar},\eps_{\hbar})^{J+1} }a\circ G_{\eps}^{t\tau_{\hbar}}(\rho)d\eps\right\},$$
and
$$A_+(\hbar,t,\delta):=\sup_{E\in[E_1-\delta/2,E_2+\delta/2]}\sup_{\rho\in T_E^*M}\left\{\frac{1}{(2\eps_{\hbar})^{J+1}} \int_{(-\eps_{\hbar},\eps_{\hbar})^{J+1} }a\circ G_{\eps}^{t\tau_{\hbar}}(\rho)d\eps\right\},$$
 where $G_{\eps}^t$ is the Hamiltonian flow associated to the Hamiltonian $p_{\eps}(x,\xi)=p_0(x,\xi)+ V(\eps,x)$. The following result is the key step to connect the 
quantum quantities to their classical counterparts:
\begin{prop}\label{p:semiclassical-reduction} Define 
$$\Lambda^{1/2}_{\text{max}}:=\limsup_{t\rightarrow+\infty}\frac{1}{|t|}\log\left(\sup_{\rho\in T_{1/2}^*M}\|d_{\rho}G_0^t\|\right).$$ 
Let $(\psi_{\hbar})_{0<\hbar\leq 1}$ be a sequence of normalized states satisfying~\eqref{e:energy-loc-states}, and let $\eps_{\hbar}\rightarrow 0$. Take 
$\nu_0 \in [0, 1/2)$ and define
$$\tau_{\hbar}:= \frac{\nu_0  |\log\hbar|}{\sqrt{2E_2}\Lambda^{1/2}_{\max}}.$$

Then, there exists $\delta_1>0$ such that for every $a$ in $\ml{C}^{\infty}_c(T^*M)$, for every $t\in[0,1]$, and for every $0<\delta\leq\delta_1$, one has, as $\hbar\rightarrow 0^+$ that
\begin{equation}
 A_-(\hbar,t,\delta)+o(1)\leq\frac{1}{(2\eps_{\hbar})^{J+1}} \int_{(-\eps_{\hbar},\eps_{\hbar})^{J+1} }\left\la u_{\hbar}^{\eps}(t\tau_{\hbar}),\Oph(a)
u_{\hbar}^{\eps}(t\tau_{\hbar})\right\ra d\eps\leq A_+(\hbar,t,\delta)+o(1),
\end{equation}
where the remainder is uniform for $t\in[0,1]$, and $u_{\hbar}^{\eps}(t')$ is the solution at time $t'$ of~\eqref{e:pert-schrodinger} with initial data $\psi_{\hbar}$.

\end{prop}

\begin{rema}
We underline that this proposition is valid for any compact Riemannian manifold (without any restriction on their sectional curvature). In the case where $M$ is a 
negatively curved surface, one has $\Lambda_{\max}^{1/2}\leq U_+$. Moreover, if the sectional curvature is constant with $K\equiv -1$, one has $\Lambda_{\max}^{1/2}=U_+=1$.
\end{rema}

Using this proposition, one can see that proving Theorem~\ref{t:maintheo} is reduced to the computation of 
$\lim_{\hbar\rightarrow 0}A_-(\hbar,t,\delta)$ and of $\lim_{\hbar\rightarrow 0}A_-(\hbar,t,\delta)$. This question relies only on aspects related to the 
classical dynamics of the geodesic flow, and it will be the main purpose of sections~\ref{s:struct-stab} and~\ref{s:mixing}. The proof of theorem~\ref{t:maintheo} is in fact the direct combination of the previous proposition with corollary~\ref{c:maincoro}.

\begin{proof}

The proof of this proposition follows standard arguments from semiclassical analysis. We introduce $0\leq \chi_{\delta}\leq 1$ a smooth function on 
$\IR$ which is equal to $1$ on the interval $[E_1,E_2]$ and $0$ outside some interval $[E_1-\delta/2,E_2+\delta/2]$ (say for some $\delta<E_1/6$). Thanks to the energy 
localization of $\psi_{\hbar}$, we can then write
$$\psi_{\hbar}=\chi_{\delta}(\hat{P}_0(\hbar))\psi_{\hbar}+o(1).$$
Recall that the operator $\chi_{\delta}(\hat{P}_{\eps}(\hbar))$ is (modulo $\ml{O}(\hbar^{\infty})$) a pseudodifferential operator with principal 
symbol $\chi_{\delta}\circ p_{\eps}(x,\xi)$. Then, we find that, uniformly for $\eps\in(-\eps_{\hbar},\eps_{\hbar})^{J+1}$,
$$\psi_{\hbar}=\chi_{\delta}(\hat{P}_{\eps}(\hbar))\psi_{\hbar}+o(1)+\ml{O}(\eps_{\hbar}).$$
This implies that, one has, uniformly for $t\in[0,1]$,
$$\mu_{\hbar}^{\eps_{\hbar}}(t)(a)=\frac{1}{(2\eps_{\hbar})^{J+1}} \int_{(-\eps_{\hbar},\eps_{\hbar})^{J+1} }
\left\la e^{+\frac{it\tau_{\hbar}\hat{P}_{\eps}(\hbar)}{\hbar}}\Oph(a)\chi_{\delta}(\hat{P}_{\eps}(\hbar))e^{-\frac{it\tau_{\hbar}\hat{P}_{\eps}(\hbar)}{\hbar}}\psi_{\hbar}\right\ra d\eps+o(1).$$
We can now apply the composition formula for pseudodifferential operators and we find that
$$\mu_{\hbar}^{\eps_{\hbar}}(t)(a)=\frac{1}{(2\eps_{\hbar})^{J+1}} \int_{(-\eps_{\hbar},\eps_{\hbar})^{J+1} }
\left\la e^{+\frac{it\tau_{\hbar}\hat{P}_{\eps}(\hbar)}{\hbar}}\Oph(a\times\chi_{\delta}\circ p_{\eps})e^{-\frac{it\tau_{\hbar}\hat{P}_{\eps}(\hbar)}{\hbar}}\psi_{\hbar}\right\ra d\eps+o(1),$$
where the remainder is still uniform in $t$. For $\hbar$ small enough and $\eps\in(-\eps_{\hbar},\eps_{\hbar})^{J+1}$, the function $a\times\chi_{\delta}\circ p_{\eps}$ is compactly 
supported in the energy layer $\{(x,\xi):E_1-\delta\leq|\xi|^2/2\leq E_2+\delta\}$. In particular, if we take $\delta>0$ small enough in a way that depends only 
on $\nu_0  \in [0,1/2)$, $E_1$ and $E_2$, we can apply Egorov Theorem up to the time $\tau_{\hbar}$ (uniformly for $\eps\in(-\eps_{\hbar},\eps_{\hbar})^{J+1}$) -- see 
Theorem~\ref{t:large-time-egorov3} from the appendix. More precisely, we have that, for $\delta>0$ small enough, 
$$\mu_{\hbar}^{\eps_{\hbar}}(t)(a)=\frac{1}{(2\eps_{\hbar})^{J+1}} \int_{(-\eps_{\hbar},\eps_{\hbar})^{J+1} }\left\la\psi_{\hbar},
\Oph(a\circ G_{\eps}^{t\tau_{\hbar}}\times\chi_{\delta}\circ p_{\eps})\psi_{\hbar}\right\ra d\eps+o(1),$$
where the remainder is still uniform for $t$ in $[0,1]$. As discussed in appendix~\ref{a:pdo}, the symbol 
$a\circ G_{\eps}^{t\tau_{\hbar}}\times\chi_{\delta}\circ p_{\eps}$ remains in a class of symbols $S^{-\infty,0}_{\nu_0'}(T^*M)$ with $0\leq \nu_0<\nu_0'<1/2$ 
(with semi-norms which are uniformly bounded for $t\in[0,1]$ and $\eps\in(-\eps_{\hbar},\eps_{\hbar})^{J+1}$). 
In particular, we can use the results from paragraph~\ref{ss:antiwick} and replace $\Oph$ by a 
positive quantization $\Oph^+$ (see~\eqref{e:positive-quantization} for instance), i.e.
$$\mu_{\hbar}^{\eps_{\hbar}}(t)(a)=\frac{1}{(2\eps_{\hbar})^{J+1}} \int_{(-\eps_{\hbar},\eps_{\hbar})^{J+1} }\left\la\psi_{\hbar},
\Oph^+(a\circ G_{\eps}^{t\tau_{\hbar}}\times\chi_{\delta}\circ p_{\eps})\psi_{\hbar}\right\ra d\eps+o(1),$$
where the remainder is still uniform for $t$ in $[0,1]$. We use Calder\'on-Vaillancourt Theorem one more time and the fact that 
$\|\chi_{\delta}\circ p_{\eps}-\chi_{\delta}\circ p_{0}\|_{\infty}=\ml{O}(\eps_{\hbar})$ uniformly for $\eps\in(-\eps_{\hbar},\eps_{\hbar})^{J+1}$. It allows us to write
$$\mu_{\hbar}^{\eps_{\hbar}}(t)(a)=\frac{1}{(2\eps_{\hbar})^{J+1}} \int_{(-\eps_{\hbar},\eps_{\hbar})^{J+1} }\left\la\psi_{\hbar},
\Oph^+(a\circ G_{\eps}^{t\tau_{\hbar}}\times\chi_{\delta}\circ p_{0})\psi_{\hbar}\right\ra d\eps+o(1).$$
The result follows then from the fact that
$$A_-(\hbar,t,\delta)\chi_{\delta}\circ p_{0}\leq \frac{1}{(2\eps_{\hbar})^{J+1}} \int_{(-\eps_{\hbar},\eps_{\hbar})^{J+1} } a\circ G_{\eps}^{t\tau_{\hbar}}\times\chi_{\delta}\circ p_{0} d\eps
\leq A_+(\hbar,t,\delta)\chi_{\delta}\circ p_{0}.$$
\end{proof}

\section{Strong structural stability of Anosov flows}
\label{s:struct-stab}
Fix $a \in C^{\infty}(T_{1/2}^*M)$ and denote by $\tilde{a}$ its $0$-homogeneous extension to $T^*M-M$. Thanks to 
proposition~\ref{p:semiclassical-reduction} and in order to prove theorem~\ref{t:maintheo}, we now have to understand the asymptotic behaviour of 
integrals of the following form:
$$I_{x_0,\xi_0}(b,T_0):=\frac{1}{(2b)^{J+1}}\int_{(-b,b)^{J+1}}\tilde{a}\circ G_{\eps}^{T_0}(x_0,\xi_0)d\eps,$$
when $T_0\rightarrow+\infty$, $b\rightarrow 0$ (simultaneously) and $p_0(x_0,\xi_0)\in [E_1-\delta/2,E_2+\delta/2]$ (for some small $\delta>0$). Our goal is to prove that, for an adapted 
range of $b$ and $T_0$, this quantity converges to $\int_{T^*_{1/2}M} adL$ \emph{uniformly} for $(x_0,\xi_0)$ satisfying $p_0(x_0,\xi_0)\in [E_1-\delta/2,E_2+\delta/2]$.

Our strategy to study this integral will be to conjugate the perturbed flow $G_{\eps}^{t}$ to the unperturbed one $G_0^{t}$ and to use the ergodic properties of $G_0^t$. Recall that two Anosov flows which are ``close'' are conjugated up to a reparametrization of time~\cite{Ano67,dLLMM86, KaHa}. The homeomorphism that conjugates the two flows has a 
priori low regularity; however, we will take advantage of the fact that this homeomorphism depends in a smooth way of the perturbation parameter $\eps$. To our 
knowledge, the smooth dependence in the perturbation parameter was observed by de la Llave, Marco and Moriyon in~\cite{dLLMM86} (appendix $A$) building on earlier 
proofs of Moser and Mather for diffeomorphisms~\cite{Mo69, Sm67}. Proving this regularity relies on 
the use of an implicit function theorem on some manifolds of mappings. Following their approach, the purpose of this section is to describe some geometric 
properties of the conjugating homeomorphism that they obtain.

The main result of this section is proposition~\ref{p:smooth-curve} which describes 
in some sense the transversality of the conjugating homeomorphism with respect to the weakly stable foliation. This
 propositions will allow us to reduce the problem of the convergence 
of the integrals $I_{x_0,\xi_0}(b,T_0)$ into a problem of equidistribution of small pieces of unstable manifolds under the geodesic flow $G_0^t$; this reduction will be the object of section~\ref{s:mixing}.\\

We now fix $(x_0,\xi_0)$ in $T_E^*M$, where $E_1-\delta/2\leq E=p_0(x_0,\xi_0)\leq E_2+\delta/2$ (with $\delta>0$ small enough).  Observe that, for a fixed choice of $(x_0,\xi_0)$ in $T^*M$ and for every $\eps$, the trajectories $G_{\eps}^{T_0}(x_0,\xi_0)$ remains a priori in different 
$3$~dimensional submanifolds inside $T^*M$. Our first step will be to project everything on the unit cotangent bundle $T_{1/2}^*M$ -- see paragraph~\ref{ss:newflow}. 
This will define a new flow on $T_{1/2}^*M$, close to that of $X_0$, to which we will apply the strong structural stability theorem.

\subsection{A new flow on $T_{1/2}^*M$} \label{ss:newflow}

Our first step is to define a new vector field on $T_{1/2}^*M$ which is ``close'' to the geodesic vector field $X_0$ and which allows to rewrite the integral 
$I_{x_0,\xi_0}(b,T_0)$ into an integral involving only quantities acting on $T_{1/2}^*M$. For that purpose, we introduce
$$\Sigma_{x_0,\xi_0}^{\eps}:=\{(x,\xi)\in T^*M:p_{\eps}(x,\xi)=p_{\eps}(x_0,\xi_0)\},$$ 
which is an energy layer for the Hamiltonian $p_{\eps}$. For every $t$ in $\IR$, $G_{\eps}^t(x_0,\xi_0)$ remains inside the energy layer $\Sigma_{x_0,\xi_0}^{\eps}$. 
We can now introduce two diffeomorphisms:
$$\theta^{\eps}_{x_0,\xi_0}:\Sigma_{x_0,\xi_0}^{\eps}\rightarrow T_{1/2}^*M,\ (x,\xi)\mapsto (x,\xi/\|\xi\|),$$
and its inverse
$$\left(\theta^{\eps}_{x_0,\xi_0}\right)^{-1}:T_{1/2}^*M\rightarrow \Sigma_{x_0,\xi_0}^{\eps},\ (x,\xi)\mapsto \left(x,\sqrt{2(p_{\eps}(x_0,\xi_0)-V_{\eps}(x))}\xi\right),$$
where we have denoted by $V_{\eps}$ the function $ V(\eps,x)$.

\begin{rema}
The map $\Pi(x,\xi):=(x,\xi/\|\xi\|)$ is well-defined on $T^*M-0$. The map $\theta^{\eps}_{x_0,\xi_0}$ is the restriction of $\Pi$ to the energy layer 
$\Sigma_{x_0,\xi_0}^{\eps}$. Observe also that
$$\left(\theta^{\eps}_{x_0,\xi_0}\right)^{-1}(x,\xi):=\left(x,\sqrt{\|\xi_0\|^2+2(\eps_0(V_0(x_0)-V_0(x))+\ldots+\eps_{J}(V_J(x_0)-V_J(x)))}\xi\right).$$

\end{rema}

Thanks to these two families of maps, we can introduce a family of flows on $T_{1/2}^*M$, i.e.
$$\varphi_{\eps,x_0,\xi_0}^t=\theta^{\eps}_{x_0,\xi_0}\circ G_{\eps}^{t/\sqrt{2E}}\circ \left(\theta^{\eps}_{x_0,\xi_0}\right)^{-1},$$
where $E=p_0(x_0,\xi_0)$. Then, we can rewrite the integral $I_{x_0,\xi_0}(b,T_0)$ as follows
\begin{equation}\label{e:integral-new-flow}
 I_{x_0,\xi_0}(b,T_0)=\frac{1}{(2b)^{J+1}}\int_{(-b,b)^{J+1}}a\circ\varphi_{\eps,x_0,\xi_0}^{T_0\sqrt{2E}}\left(x_0,\frac{\xi_0}{\|\xi_0\|}\right)d\eps.
\end{equation}
The end of this section will be devoted to the application of the strong structural stability theorem to the flow $(\varphi_{\eps,x_0,\xi_0}^t)_t$. 
For that purpose, we need to compute the vector field associated to this new flow:
\begin{equation}
Y_{x_0,\xi_0}^{\eps}(\rho):=\frac{d}{dt}\left(\varphi_{\eps,x_0,\xi_0}^t(\rho)\right)_{t=0}. 
\end{equation}
Using the notations of section~\ref{s:geom-back}, the following holds:
\begin{lemm} One has, for every $\rho=(x,\xi)$ in $T_{1/2}^*M$,
\begin{equation}\label{e:newflow}
Y_{x_0,\xi_0}^{\eps}(\rho)
=c_{\eps,x_0,\xi_0}(x)X_0(\rho)+\frac{1}{\sqrt{2E}c_{\eps,x_0,\xi_0}(x)}\sum_{j=0}^J\eps_jg_x^*\left(d_xV_j,\xi^{\perp}\right)J_{\rho}W_{1/2}(\rho),
\end{equation}
 where $c_{\eps,x_0,\xi_0}(x):=\sqrt{\frac{p_{\eps}(x_0,\xi_0)-V_{\eps}(x)}{p_0(x_0,\xi_0)}}$, and $E:=p_0(x_0,\xi_0)$.
\end{lemm}

\begin{rema} 
We can reestablish the dependence in $(x_0,\xi_0)$ more clearly and write:
$$c_{\eps,x_0,\xi_0}(x)=\sqrt{1+\frac{2}{\|\xi_0\|^2}\sum_{j=0}^J\eps_j(V_j(x_0)-V_j(x))}=1+\ml{O}_{x,x_0,\xi_0}(\|\eps\|).$$
\end{rema}

\begin{proof} By definition, one has
$$Y^{\eps}_{x_0,\xi_0}:=\frac{d\theta^{\eps}_{x_0,\xi_0}.\left(X_{\eps}\circ\left(\theta^{\eps}_{x_0,\xi_0}\right)^{-1}\right)}{\sqrt{2E}},$$
where
$$X_{\eps}=X_0+\eps_0 X_{V_0}+ \eps_1X_{V_1}+\ldots+\eps_JX_{V_J},$$
with $X_0$ is geodesic vector field and $X_{V_j}$ the Hamiltonian vector field associated to $V_j$. In particular, using the results of section~\ref{s:geom-back}, we have
$$\frac{d\theta^{\eps}_{x_0,\xi_0}.
X_{0}\circ\left(\theta^{\eps}_{x_0,\xi_0}\right)^{-1}}{\sqrt{2E}}=\sqrt{\frac{p_{\eps}(x_0,\xi_0)-V_{\eps}}{E}}X_0.$$
We will now express
$$\frac{d\theta^{\eps}_{x_0,\xi_0}.
X_{V_j}\circ\left(\theta^{\eps}_{x_0,\xi_0}\right)^{-1}}{\sqrt{2E}}.$$
For that purpose, we recall from section~\ref{s:geom-back} that $X_{V_j}\circ\left(\theta^{\eps}_{x_0,\xi_0}\right)^{-1}$ belongs to the vertical bundle $\ml{V}:=\cup_{\rho}\ml{V}_{\rho}$. The action of 
the tangent map $d_{\rho}\Pi$ on the vertical space is given by $Y_0(\rho)\mapsto 0$ and $J_{\rho}W_{\rho}\mapsto J_{\Pi(\rho)}W_{\Pi(\rho)}$. So, for 
every $\rho=(x,\xi)$ in $T_{1/2}^*M$, the vector we are interested in can be expressed as follows:
$$\frac{d\theta^{\eps}_{x_0,\xi_0}.
X_{V_j}\circ\left(\theta^{\eps}_{x_0,\xi_0}\right)^{-1}}{\sqrt{2E}}=\frac{1}{\sqrt{p_{\eps}(x_0,\xi_0)-V_{\eps}(x)}}g_x^*\left(d_xV_j,\xi^{\perp}\right)J_{\rho}W_{1/2}(\rho).$$
Finally, we can write, for every $\rho=(x,\xi)$ in $T_{1/2}^*M$,
$$Y_{x_0,\xi_0}^{\eps}(\rho)
=c_{\eps,x_0,\xi_0}(x)X_0(\rho)+\frac{1}{\sqrt{2E}c_{\eps,x_0,\xi_0}(x)}\sum_{j=0}^J\eps_jg_x^*\left(d_xV_j,\xi^{\perp}\right)J_{\rho}W_{1/2}(\rho).$$
\end{proof}


We denote by $\ml{V}^2(T_{1/2}^*M)$ the space of $\ml{C}^2$ vector fields on $T_{1/2}^*M$. For every open neighborhood $\ml{U}$ of $X_0$ in $\ml{V}^2(T_{1/2}^*M)$, one can verify 
that there exists $b(\ml{U})>0$ such that, for every $0<b<b(\ml{U})$ and for every $(x_0,\xi_0)$ satisfying $p_0(x_0,\xi_0)\in[E_1-\delta/2,
E_2+\delta/2]$, $Y_{x_0,\xi_0}^{\eps}$ belongs to $\ml{U}$. One can also verify that the map
$$\eps\in(-b,b)^{J+1}\mapsto Y^{\eps}_{x_0,\xi_0}\in\ml{V}^2(T_{1/2}^*M)$$
is of class $\ml{C}^{\infty}$ with respect to the parameter $\eps$. 


\subsection{Strong structural stability of Anosov vector fields}

The previous paragraph has allowed us to project on $T_{1/2}^*M$ the trajectories of the perturbed flow $G_{\eps}^t$ thanks to the family of vector fields 
$Y^{\eps}_{x_0,\xi_0}$ (which are small perturbations of $X_0$). We will now recall strong structural stability property for Anosov flows in the same way as it was stated in~\cite{dLLMM86} -- see also~\cite{KKPW89}. For that purpose, 
we introduce manifolds of mappings that will be involved in this theorem -- see appendix~\ref{a:map-manifold} for a brief reminder on 
their differential structure.

We introduce the following space
$$\mathcal{C}_{X_0}(T_{1/2}^*M):=\left\{h\in\ml{C}^{0}(T_{1/2}^*M,T_{1/2}^*M):\ \forall \rho\in T_{1/2}^*M,
\ \left(\frac{d}{dt}h\circ G_0^t(\rho)\right)_{t=0}=D_{X_0}h\ \text{exists}\right\},$$
which can be endowed with a smooth differential structure modeled on the Banach spaces of continuous sections $s:T_{1/2}^*M\mapsto h^*TT_{1/2}^*M$ which are differentiable 
along the geodesic flow. This manifold contains an ``adapted'' submanifold $\ml{M}$ which contains $\text{Id}_{T_{1/2}^*M}$ in its interior and \emph{for which the elements 
are in some sense ``transversal'' to the geodesic vector field} $X_0$ -- see appendix $A$ of~\cite{dLLMM86} or appendix~\ref{a:map-manifold} of the present article. 

The structural stability theorem can be then stated as follows (Theorem $A.2$ in~\cite{dLLMM86}):

\begin{theo}\label{t:stability}[Strong structural stability] Assume $X_0$ is an Anosov vector field. There exists an open neighborhood $\ml{U}_0(X_0)$ of $X_0$ in $\ml{V}^2(T_{1/2}^*M)$ and a 
unique $\ml{C}^2$ map $S_0:\ml{U}_0(X_0)\rightarrow \ml{M}\times\ml{C}^0(T_{1/2}^*M,\IR)$ such that $S_0(X_0)=(\text{Id}_{T_{1/2}^*M},1)$ and if $S_0(X)=(h,\tau)$, then
\begin{equation}
 \label{e:implicit-equation}
D_{X_0}h-\tau X\circ h=0_{T_{1/2}^*M}(h),
\end{equation}
where $0_{T_{1/2}^*M}$ is the zero section.
\end{theo}

\begin{rema}\label{r:reparametrization} In fact, the neighborhood $\ml{U}_0(X_0)$ can be chosen small enough to ensure that $h$ 
is an homeomorphism -- appendix $A$ of~\cite{dLLMM86}. Surjectivity follows from the facts that the identity map is of topological degree $1$ and that the degree is invariant by homotopy. 
As $T_{1/2}^*M$ is compact and as $h$ is continuous, it remains to prove the injectivity. For that purpose, one can observe that, if we denote $G_X^t$ the flow 
associated to the vector field $X$, then for every $t$ in $\IR$, the following holds:
$$h\circ G_0^t(\rho)=G_X^{\tilde{\tau}_{\rho}(t)}\circ h(\rho),$$
where $\tilde{\tau}_t(\rho)=\int_0^t\tau\circ G_0^s(\rho)ds$. Then, one can use expansivity of Anosov flows. Precisely, one can find $\eta_0>0$ 
small such that every $\rho$ in $T_{1/2}^*M$ satisfies the following: if for $\rho'\in T_{1/2}^*M$, there exists $\tau_0:\IR\rightarrow \IR$ a continuous surjection with 
$\tau_0(0)=0$ such that
$$\forall t\in\IR,\ d(G_0^t(\rho),G_0^{\tau_0(t)}(\rho'))\leq\eta_0,$$
there exists $t_0$ such that $G_0^{t_0}(\rho)=\rho'$~\cite{Rug07}. Combining this property to the above reparametrization formula allows to prove injectivity 
(at least for a small enough neighborhood of $X_0$). Using that $h$ is an homeomorphism, we can also write the following formula connecting the two flows:
\begin{equation}\label{e:struc-stab-useful} 
\forall\ t\in\IR,\ h\circ G_0^{\tau(t,\rho)}\circ h^{-1}(\rho)=G_X^t(\rho),
\end{equation}
where
$$\tau(t,\rho):=\int_{0}^t\frac{ds}{\tau\circ h^{-1}\circ G_X^s(\rho)}.$$
In section~\ref{s:mixing}, we will use strong structural stability under this form.
\end{rema}

From the construction of the vector field $Y_{x_0,\xi_0}^{\eps}$, we know that there exists some $b_0>0$ such that, for every $\eps\in(-b_0,b_0)^{J+1}$ and for every 
$(x_0,\xi_0)$ satisfying $E_1-\delta/2\leq p_0(x_0,\xi_0)\leq E_2+\delta/2$, $Y_{x_0,\xi_0}^{\eps}$ belongs to the neighborhood $\ml{U}_0(X_0)$ of the theorem. We define 
then a family of curves on $T_{1/2}^*M$ paramatrized by $\eps\in(-b_0,b_0)^{J+1}$, namely
\begin{equation}\label{e:family-curve}
\alpha_{x_0,\xi_0}(\eps)=\left(h_{x_0,\xi_0}^{\eps}\right)^{-1}\left(x_0,\frac{\xi_0}{\|\xi_0\|}\right),
\end{equation}
where $(h_{x_0,\xi_0}^{\eps},\tau_{x_0,\xi_0}^{\eps})=S_0(Y_{x_0,\xi_0}^{\eps})$ and $(x_0,\xi_0)$ varies in the allowed energy layers. We can now 
rewrite the integral $I_{x_0,\xi_0}(b,T_0)$ using~\eqref{e:integral-new-flow},~\eqref{e:struc-stab-useful} and~\eqref{e:family-curve}. We obtain
\begin{equation}\label{e:integral-rewritten}
 I_{x_0,\xi_0}(b,T_0)=\frac{1}{(2b)^{J+1}}\int_{(-b,b)^{J+1}} a\circ h_{x_0,\xi_0}^{\eps}\circ G_0^{\tau_{x_0,\xi_0}^{\eps}(T_0\sqrt{2E},x_0,\xi_0/\|\xi_0\|)} 
\left(\alpha_{x_0,\xi_0}(\eps)\right)d\eps,
\end{equation}
where
$$\tau_{x_0,\xi_0}^{\eps}\left(T_0\sqrt{2E},x_0,\frac{\xi_0}{\|\xi_0\|}\right):=\int_{0}^{T_0\sqrt{2E}}\frac{ds}{\tau_{x_0,\xi_0}^{\eps}\circ (h_{x_0,\xi_0}^{\eps})^{-1}
\circ \varphi_{\eps,x_0,\xi_0}^s(x_0,\xi_0/\|\xi_0\|)}.$$
Thus, in order to study the asymptotic behaviour of $I_{x_0,\xi_0}(b,T_0)$, we will need to understand the action of $G_0^t$ on the subsets defined by the maps 
$\eps\mapsto\alpha_{x_0,\xi_0}(\eps)$. For that purpose, we need to understand the geometric properties of these curves. More specifically, we need to quantify 
their transversality to the weakly stable foliation. This geometric description is the main objective of this section.

\begin{rema}
 In section~\ref{s:mixing}, we will also have to check that several constants are uniformly bounded in terms of $(x_0,\xi_0)$, in the allowed energy layer. A simple 
observation that we can also make is that there exists some uniform constant $C>0$ depending only on the size of the neighborhood $\ml{U}_0(X_0)$ such that
\begin{equation}\label{e:taylor-stability-equation}
 \max\left\{d_{\ml{C}^0}\left(h_{x_0,\xi_0}^{\eps},\text{Id}\right),d_{\ml{C}^0}\left((h_{x_0,\xi_0}^{\eps})^{-1},\text{Id}\right),\left\|\tau_{x_0,\xi_0}^{\eps}-1\right\|_{\ml{C}^0}\right\}\leq 
C\left\|Y_{x_0,\xi_0}^{\eps}-X_0\right\|_{\ml{C}^2}.
\end{equation}
Thanks to equation~\eqref{e:newflow}, we can observe that the right-hand side of this inequality is $\ml{O}(\|\eps\|)$ where the involved constant is uniform for $(x_0,\xi_0)$ 
in the allowed interval of energy.
\end{rema}

\begin{rema}\label{r:inv-not-smooth} At this point, we would like to make an important observation. In section~\ref{s:mixing}, one of the difficulty we will encounter is that 
the map $\eps\mapsto\alpha_{x_0,\xi_0}(\eps)$ is not a priori of class $\ml{C}^1$. In fact, the strong structural stability theorem stated above only implies that $\eps\mapsto h^{\eps}_{x_0,\xi_0}$ is of class $\ml{C}^1$ 
as a map from a neighborhood of $0$ in $\mathbb{R}^{J+1}$ to $\ml{C}^0(T_{1/2}^*M,T_{1/2}^*M)$ but there is a priori no reason that the map $\eps\mapsto (h^{\eps}_{x_0,\xi_0})^{-1}$ would be of class 
$\ml{C}^1$. In some sense, it is related to the fact that $h\mapsto h^{-1}$ is only $\ml{C}^0$ as a map from $\ml{C}^0(T_{1/2}^*M,T_{1/2}^*M)$ to itself while it 
would be of class $\ml{C}^1$ as a map from $\ml{C}^1(T_{1/2}^*M,T_{1/2}^*M)$ to $\ml{C}^0(T_{1/2}^*M,T_{1/2}^*M)$. We cannot expect $ h^{\eps}_{x_0,\xi_0}$ to be 
a $\ml{C}^1$ map as the conjugating homeomorphism in the strong structural stability theorem is in general only of H\"older class. On the other hand, one could have 
expected that revisiting the proof of the strong structural stability would provide an information on the regularity of the map $X\mapsto h^{-1}$ but we did not find any 
way of proving such statement. In fact, it even seems to us that one cannot expect any smoothness dependence of 
the map $X\mapsto h^{-1}$ under such generality\footnote{This problem is not specific to flows and the same difficulty seems to occur in the case of 
strong structural stability for diffeomorphisms.}.
\end{rema}

In order to solve the problem raised in remark~\ref{r:inv-not-smooth}, we will approximate $\alpha_{x_0,\xi_0}$ by a ``good'' map $\eps\mapsto\tilde{\alpha}_{x_0,\xi_0}(\eps)$ which is of class $\ml{C}^{1}$. For that purpose, we introduce the following vector field on $S^*M$:
\begin{equation}\label{e:derivative-heps} Z_{\eps}(\rho):=-\sum_{j=0}^J\eps_j\left(\ml{L}^s(b_1^j) X^s_{1/2}+\ml{L}^u(b_1^j)X^u_{1/2}\right)(\rho).\end{equation}
where, for every $0\leq j\leq J$, we set 
$$b_1^j(x,\xi):=\frac{g_x^*(d_xV_j,\xi^{\perp})}{\|\xi_0\|},$$
and, for every continuous function $b$ on $S^*M$, 
\begin{equation}\label{e:averaging-unst}\ml{L}^u(b):=\int_{0}^{+\infty}\left(\frac{b}{U^u-U^s}\right)\circ G_0^{t} e^{-\int_0^{t}U^u\circ G_0^sds}dt,\end{equation}
and
\begin{equation}\label{e:averaging-st}\ml{L}^s(b):=\int_{0}^{+\infty}\left(\frac{b}{U^u-U^s}\right)\circ G_0^{-t} e^{-\int_0^{-t}U^s\circ G_0^sds}dt.\end{equation}
Note that $Z_{\eps}$ implicitely depends on $(x_0,\xi_0)$. We recognize the quantities appearing in the definition of the admissibility operator~\eqref{e:admissibility-operator}. Using this vector field $Z_{\eps}$, we define the following map:
\begin{equation}\label{e:explicit-curve}
\tilde{\alpha}_{x_0,\xi_0}(\eps):=\exp_{x_0,\frac{\xi_0}{\|\xi_0\|}}\left(Z_{\eps}\right), 
\end{equation}
where $\exp$ is the exponential map induced by the Riemannian structure on $T_{1/2}^*M$. The following proposition is the main result of this section:

\begin{prop}\label{p:smooth-curve} Suppose $M$ is a smooth compact Riemannian surface of negative curvature.
Let $0<\gamma_1<\gamma_c:=\frac{U_-}{U_++U_-}$. For every $(x_0,\xi_0)$ satisfying $p_0(x_0,\xi_0)\in[E_1-\delta/2,E_2+\delta/2]$ and for every $\eps\in(-b_0,b_0)^{J+1}$, one has
$$d_{T_{1/2}^*M}(\alpha_{x_0,\xi_0}(\eps),\tilde{\alpha}_{x_0,\xi_0}(\eps))\leq C\|\eps\|^{1+\gamma_1}.$$
Moreover, $b_0$ and $C$ can be chosen uniformly for $(x_0,\xi_0)$ in the allowed energy layers.
\end{prop}

We observe that, in constant curvature $K\equiv-1$, one has $\gamma_c=\frac{1}{2}$. We also underline that our admissibility assumption~\eqref{e:admissibility} on the 
family of potentials exactly says that at least one the partial derivative is nonzero. In other words, \emph{we get some transversality to the weakly stable foliation 
for a large choice of $\eps$ in the box $(-b_0,b_0)^{J+1}$}. 



\subsection{Proof of proposition~\ref{p:smooth-curve}}

The proof of this proposition is somewhat technical as we have to deal with the machinery of 
Banach manifolds. Yet, the main lines of the proof can be easily explained. First, the fact that $\eps\mapsto h_{x_0,\xi_0}^{\eps}$ is of class $\ml{C}^{2}$ means that 
its ``linearization'' $\eps\mapsto v_{x_0,\xi_0}^{\eps}:=\exp^{-1} h_{x_0,\xi_0}^{\eps}$ is of class $\ml{C}^2$, where $\exp$ is the exponential map induced by the 
Riemannian metric on $T_{1/2}^*M$. We can then exploit the fact that $h_{x_0,\xi_0}^{\eps}$ 
satisfies the implicit equation~\eqref{e:implicit-equation} in order to obtain an asymptotic expansion of $v_{\eps}$ in terms of $\eps$. 
This first step of the proof is given in paragraphs~\ref{sss:step1} and~\ref{sss:step2}. Once we have this expansion, we can denote by $-\tilde{Z}_{x_0,\xi_0}^{\eps}$ the first terms in 
the expansion and we can verify that $\exp \tilde{Z}_{x_0,\xi_0}^{\eps}$ is a good approximation of $(h_{x_0,\xi_0}^{\eps})^{-1}$ in the sense of proposition~\ref{p:smooth-curve}. This last step 
is carried out in paragraphs~\ref{sss:step3} and~\ref{sss:step4}.

\subsubsection{Differentiation of the implicit equation~\eqref{e:implicit-equation}}
\label{sss:step1}

Define the following map:
$$F:\ml{V}^2(T_{1/2}^*M)\times\ml{M}\times\ml{C}^0(T_{1/2}^*M,\IR)\rightarrow\ml{C}^0(T_{1/2}^*M,TT_{1/2}^*M),\ (X,h,\tau)\mapsto D_{X_0}h-\tau X\circ h.$$
According to~\cite{dLLMM86} (appendix $A$), this map is of class $\ml{C}^2$. From the structural stability equation, we have that, for $\eps$ small enough,
\begin{equation}\label{e:implicit}
F(Y^{\eps}_{x_0,\xi_0},h^{\eps}_{x_0,\xi_0},\tau^{\eps}_{x_0,\xi_0})=0_{T_{1/2}^*M}(h^{\eps}_{x_0,\xi_0}).
\end{equation}

\begin{rema}
In order to alleviate notations, we will often omit the subscript $(x_0,\xi_0)$ in this paragraph and restablish it everytime it will be necessary.
\end{rema}

We introduce the map $0: h\in\ml{M}\mapsto 0_{T_{1/2}^*M}\circ h\in \ml{C}^0(T_{1/2}^*M,TT_{1/2}^*M)$ which is of class $\ml{C}^{\infty}$ -- see lemma 
$A.1$ in~\cite{dLLMM86} or~\cite{Ab63,dLLO99}. This map can also be identified with a smooth map acting on $\ml{V}^2(T_{1/2}^*M)\times\ml{M}\times\ml{C}^0(T_{1/2}^*M,\IR)$.

\begin{rema}\label{r:targetspace} As explained in appendix~\ref{a:map-manifold}, the tangent space to $\text{Id}$ in $\ml{M}$ can be identified with the space $\Gamma_{X_0}(\ml{N})$ 
of continuous sections which are differentiable along the geodesic flow and which take values in the normal bundle defined in section~\ref{s:geom-back}. 
Before starting our computation, we recall from~\cite{dLLMM86} that 
\begin{equation}\label{e:linearize}
(D_{2}F-D_2 0,D_{3}F)(X_0,\text{Id},1):\Gamma_{X_0}(\ml{N})\times\ml{C}^0(T_{1/2}^*M,\IR)\rightarrow \Gamma^0(0^*TTT_{1/2}^*M)
\end{equation} 
is a linear isomorphism onto a subspace of $\Gamma^0(0^*TTT_{1/2}^*M)$ which can be identified with 
$\ml{V}^0(T_{1/2}^*M)$, where $\ml{V}^0(T_{1/2}^*M)$ is the Banach space of continuous vector fields on $T_{1/2}^*M$.

This latter identification can be made thanks to the following observations. One has $0^*TTT_{1/2}^*M=TTT_{1/2}^*M\rvert_{T_{1/2}^*M}$, and 
$T_{\rho}T_{1/2}^*M$ can be identified with a subspace of $T_{0(\rho)}TT_{1/2}^*M$ in two different ways
\begin{itemize}
 \item using the tangent map of $0_{T_{1/2}^*M}:T_{1/2}^*M\rightarrow TT_{1/2}^*M$;
 \item using the tangent map at $0$ of the inclusion map $i_{\rho}:T_{\rho}T_{1/2}^*M\rightarrow TT_{1/2}^*M$. 
\end{itemize}
Thus, $\Gamma^0(0^*TTT_{1/2}^*M)$ can be identified with $\ml{V}^0(T_{1/2}^*M)\oplus\ml{V}^0(T_{1/2}^*M)$, and the image of the linear map~\eqref{e:linearize} is exactly 
$\{0\}\oplus\ml{V}^0(T_{1/2}^*M)$ -- see lemma $A.7$ and its proof in~\cite{dLLMM86} (see also appendix~\ref{a:map-manifold} for a brief reminder).
\end{rema}

For every $0\leq j\leq J$, we start by differentiating~\eqref{e:implicit}
$$(D_{2}F-D_20,D_3F)(Y^{\eps},h^{\eps},\tau^{\eps}).\left(\frac{\partial h^{\eps}}{\partial\eps_j},\frac{\partial\tau^{\eps}}{\partial\eps_j}\right)=-D_1F(Y^{\eps},h^{\eps},\tau^{\eps}).\frac{\partial Y^{\eps}}{\partial\eps_j},$$
In order to alleviate the notations, we denote by $\ml{A}(Y,h,\tau)$ the linear map \break $(D_{2}F-D_2 0,D_3F)(Y,h,\tau)$
acting on the tangent space $T_{h}\ml{M}\times\ml{C}^0(T_{1/2}^*M,\IR)$. Thanks to the above remark (see also appendix~\ref{a:map-manifold}), we know that it is 
a continuous linear isomorphism onto $\ml{V}^{0}(T_{1/2}^*M)$ for $(Y,h,\gamma)=(X_0,\text{Id}, 1)$. Recall now that the maps $F$ and $0$ are of class $\ml{C}^2$, that 
$\eps\mapsto Y^{\eps}$ is of class 
$\ml{C}^{\infty}$ and that the map $S_0$ is of class $\ml{C}^2$. In particular, for $\eps$ small enough (in a way that does not depend 
on $(x_0,\xi_0)$ but only on the energy layers they belong to), one can deduce that 
the map $\ml{A}(Y^{\eps},h^{\eps},\tau^{\eps})$ remains invertible. Thus, one can write, for every $0\leq j\leq J$,
\begin{equation}\label{e:diff-implicit}\left(\frac{\partial h^{\eps}}{\partial\eps_j},\frac{\partial\tau^{\eps}}{\partial\eps_j}\right)
=-\ml{A}(Y^{\eps},h^{\eps},\tau^{\eps})^{-1}D_1F(Y^{\eps},h^{\eps},\tau^{\eps}).\frac{\partial Y^{\eps}}{\partial\eps_j}.
\end{equation}
We will now approximate the right hand side of this equality in the $\ml{C}^0$ topology in order to get an approximate expression for $\frac{\partial h^{\eps}}{\partial\eps_j}$.

\subsubsection{Computation of the right hand side of~\eqref{e:diff-implicit}}
\label{sss:step2}

Using expression~\eqref{e:newflow}, we have that
$$\frac{\partial Y^{\eps}}{\partial\eps_j}=b_0^jX_0+b_1^j JW+\ml{O}(\|\eps\|),$$
where $b_0^j(x,\xi)=\frac{V_j(x_0) - V_j(x)}{\|\xi_0\|^2}$ and $b_1^j(x,\xi)=\frac{g_x^*(d_xV_j,\xi^{\perp})}{\|\xi_0\|}$ are smooth functions, and the constant 
in the remainder is uniform in $(x_0,\xi_0)$.

\begin{rema}
As $X_0$ and $JW$ are smooth vector fields, one can prove that the maps 
$$h\in\ml{C}^0(T_{1/2}^*M,T_{1/2}^*M)\mapsto X_0\circ h\in\ml{C}^0(T_{1/2}^*M,TT_{1/2}^*M),$$
and
$$h\in\ml{C}^0(T_{1/2}^*M,T_{1/2}^*M)\mapsto JW\circ h\in\ml{C}^0(T_{1/2}^*M,TT_{1/2}^*M)$$
are smooth -- lemma $A.1$ in~\cite{dLLMM86} or~\cite{Ab63, dLLO99}. We also recall that the inclusion map $$h\in\ml{C}^0_{X_0}(T_{1/2}^*M,TT_{1/2}^*M)\mapsto h\in\ml{C}^0(T_{1/2}^*M,TT_{1/2}^*M)$$
is smooth -- lemma $A.4$ in~\cite{dLLMM86}. 
\end{rema}

As the map $F$ is of class $\ml{C}^2$ and as the map 
$$\eps\in(-b_0,b_0)^{J+1}\mapsto (Y^{\eps},h^{\eps},\tau^{\eps})\in\ml{V}^2(T_{1/2}^*M)\times\ml{M}\times\ml{C}^0(T_{1/2}^*M,\IR)$$
is of class $\ml{C}^2$, one can deduce that the operator $\ml{A}(Y^{\eps},h^{\eps},\tau^{\eps})^{-1}D_1F(Y^{\eps},h^{\eps},\tau^{\eps})$ remains uniformly bounded. 
In particular, the right hand side of~\eqref{e:diff-implicit} can be rewritten, for every $0\leq j\leq J$, as
$$-\ml{A}(Y^{\eps},h^{\eps},\tau^{\eps})^{-1}D_1F(Y^{\eps},h^{\eps},\tau^{\eps})\left(b_0^jX_0+b_1^j JW\right)+\ml{O}(\|\eps\|)\in 
T_{h^{\eps}}\ml{M}\times\ml{C}^0(T_{1/2}^*M,\IR).$$

Using one more time the facts that the map $F$ is of class $\ml{C}^2$ and that the map $\eps\mapsto (Y^{\eps},h^{\eps},\tau^{\eps})$ is of class $\ml{C}^2$ 
in the appropriate spaces given by the strong structural stability theorem, we deduce that 
$\ml{A}(Y^{\eps},h^{\eps},\tau^{\eps})^{-1}D_1F(Y^{\eps},h^{\eps},\tau^{\eps}).\frac{\partial Y^{\eps}}{\partial\eps_j}$ can be approximated by 
\begin{equation}\label{e:stability-proof}\ml{A}(X_0,\text{Id},1)^{-1}D_1F(X_0,\text{Id}, 1).\left(b_0^jX_0+b_1^j JW\right),\end{equation}
in the $\ml{C}^0$ topology up to an error of order $\ml{O}(\|\eps\|)$. Again, we can observe that the constant in the remainder can be uniformly bounded in terms of 
$(x_0,\xi_0)$ in the allowed energy layers.\\

Let us give a more explicit description of the previous sum. We use the decomposition of remark~\ref{r:targetspace} and we have, for 
every $0\leq j\leq J$,
$$\ml{A}(X_0,\text{Id},1)^{-1}D_1F(X_0,\text{Id}, 1).\left(b_0^jX_0+b_1^j JW\right)=\ml{A}(X_0,\text{Id},1)^{-1}(0,-b_0^jX_0-b_1^j JW).$$
Thanks to our description of the Anosov decomposition in paragraph~\ref{ss:anosov}, one can write $JW=\frac{X^s_{1/2}-X^u_{1/2}}{U^u-U^s}$. 
Following the proof of lemma $A.7$ in~\cite{dLLMM86}, one can can write that, in the space $\Gamma_{X_0}(\ml{N})\times \ml{C}^0(T_{1/2}^*M,\IR)$,
$$\ml{A}(X_0,\text{Id},1)^{-1}(0,c_0X_0+c_u X^u_{1/2}+c_sX^s_{1/2})=$$ $$\left(\int_0^{+\infty}d_{G_0^{-t}}G_0^{t}(c_sX^s_{1/2})dt
-\int_{-\infty}^{0}d_{G_0^{-t}}G_0^{t}(c_uX^u_{1/2})dt,c_0\right),$$
for every continuous functions $c_0$, $c_u$ and $c_s$ on $T_{1/2}^*M$. As this result is an important step of our proof, we briefly recall the argument of~\cite{dLLMM86} 
in appendix~\ref{a:map-manifold}.

In our situation, it implies that
\begin{equation}\label{e:RHS}
\ml{A}(X_0,\text{Id},1)^{-1}D_1F(X_0,\text{Id}, 1).\left(b_0^jX_0+b_1^j JW\right)=\left(- \ml{L}^s(b_1^j) X^s_{1/2}-\ml{L}^u(b_1^j)X^u_{1/2},b^j_0\right).\end{equation}
where the maps $\ml{L}^{u/s}$ were defined in~\eqref{e:averaging-unst} and~\eqref{e:averaging-st}.

Combining~\eqref{e:diff-implicit},~\eqref{e:stability-proof}, and~\eqref{e:RHS}, we have shown that, up to an (uniform) error of 
order $\ml{O}(\|\eps\|)$ in the $\ml{C}^0(T_{1/2}^*M,TT_{1/2}^*M)$ topology, $\frac{\partial h^{\eps}}{\partial\eps_j}$ can be approximated, for every $0\leq j\leq J$, by the vector field
\begin{equation}\label{e:final-RHS}\rho\mapsto \left(\ml{L}^s(b_1^j) X^s_{1/2}+\ml{L}^u(b_1^j)X^u_{1/2}\right)(\rho).
\end{equation}

\subsubsection{Regularity of the function $\ml{L}^u(b)$}\label{sss:step3} Before we can conclude, 
we need to prove an intermediary lemma concerning the regularity of the map $\ml{L}^u(b)$. Precisely, we have the following lemma:
\begin{lemm}\label{l:holder} For any $b$ in $\ml{C}^{\infty}(T_{1/2}^*M,\IR)$, the functions $\ml{L}^u(b)$ and $\ml{L}^s(b)$ are $\gamma_1$-H\"older for every 
$$0\leq\gamma_1<\gamma_c=\frac{U_-}{U_++U_-}.$$
\end{lemm}

\begin{rema} This lemma is not very surprising as one knows that the conjugating homeomorphism in the strong structural stability theorem is of H\"older class for some 
small enough exponent $\gamma_1>0$. A similar idea was already exploited in~\cite{KKPW89} where the authors showed that theorem~\ref{t:stability} can be improved (modulo some small restrictions) to show that the map $S_0$ still has good regularity properties if the target space 
is the space of H\"older functions $\ml{C}^{\gamma_1}(T_{1/2}^*M,T_{1/2}^*M)$ with $\gamma_1>0$ small enough.
 
\end{rema}

\begin{proof} We only treat the case of $\ml{L}^u$. Recall that $U^u$ and $U^s$ are of class $\ml{C}^1$ due to the regularity of the 
stable/unstable foliation~\cite{PuHi75, Has94}. From 
remark~\ref{r:bound-derivative-flow}, one can observe that there exists some uniform constant $C_0>0$ such that, for every $t\geq 0$, 
$d(G_0^t(\rho),G_0^t(\rho'))\leq C_0 e^{t U_+}d(\rho,\rho')$.

We fix $0<\gamma_0<1$ and, for $\rho$ and $\rho'$ in $T_{1/2}^*M$ satisfying $d(\rho,\rho')\leq 1/2$,  define
$$T_{0}(\rho,\rho'):=\frac{1-\gamma_0}{U_+}\log\left(\frac{1}{d(\rho,\rho')}\right)\in\IR_+\cup\{+\infty\}.$$ 
We write then
$$\ml{L}^u(b)(\rho)=\int_{0}^{T_{0}(\rho,\rho')}\left(\frac{b}{U^u-U^s}\right)\circ G_0^{t}(\rho) e^{-\int_0^{t}U^u\circ G_0^s(\rho)ds}dt$$
$$+\int_{T_{0}(\rho,\rho')}^{+\infty}\left(\frac{b}{U^u-U^s}\right)\circ G_0^{t}(\rho) e^{-\int_0^{t}U^u\circ G_0^s(\rho)ds}dt.$$
We can split $\ml{L}^u(b)(\rho')$ in a similar manner. In both cases, the second term in the right hand side is bounded by $C_1e^{-U_{-}T_{0}(\rho,\rho')}$ where $C_1$ is 
some positive constant that depends only on $b$ and on the manifold. In particular, this term is bounded by $C_1 d(\rho,\rho')^{\frac{(1-\gamma_0)U_{-}}{U_+}}$. Thus, 
we have
$$\left|\ml{L}^u(b)(\rho)-\ml{L}^u(b)(\rho')\right|\leq\left|\ml{L}^u_{T_{0}(\rho,\rho')}(b)(\rho)-\ml{L}^u_{T_{0}(\rho,\rho')}(b)(\rho')\right|+2C_1 
d(\rho,\rho')^{\frac{(1-\gamma_0)U_{-}}{U_+}},$$
where
$$\ml{L}^u_{T_{0}(\rho,\rho')}(b)(\rho):=\int_{0}^{T_{0}(\rho,\rho')}\left(\frac{b}{U^u-U^s}\right)\circ G_0^{t}(\rho) e^{-\int_0^{t}U^u\circ G_0^s(\rho)ds}dt.$$
By the definition of $T_{0}(\rho,\rho')$, one can verify that
$$\left|\ml{L}^u_{T_{0}(\rho,\rho')}(b)(\rho)-\ml{L}^u_{T_{0}(\rho,\rho')}(b)(\rho')\right|\hspace{9cm}$$
$$\leq
\left|\int_{0}^{T_{0}(\rho,\rho')}\left(\frac{b}{U^u-U^s}\right)\circ G_0^{t}(\rho') e^{-\int_0^{t}U^u\circ G_0^s(\rho)ds}dt-\ml{L}^u_{T_{0}(\rho,\rho')}(b)(\rho')\right|
+C_2d(\rho,\rho')^{\gamma_0},$$
where $C_2$ depends only on $b$ and the manifold. We can remark that, for $0\leq t\leq T_{0}(\rho,\rho')$,
$$\left|\int_0^{t}U^u\circ G_0^s(\rho)ds-\int_0^{t}U^u\circ G_0^s(\rho')ds\right|\leq \frac{C_3(1-\gamma_0)}{U_+}d(\rho,\rho')^{\gamma_0}\log\left(\frac{1}{d(\rho,\rho')}\right),$$
for some uniform constant $C_3$ that depends only on the manifold. This implies that
$$e^{\int_0^{t}U^u\circ G_0^s(\rho)ds-\int_0^{t}U^u\circ G_0^s(\rho')ds}=1+\ml{O}(1)d(\rho,\rho')^{\gamma_0}\log\left(\frac{1}{d(\rho,\rho')}\right),$$
where the constant in the remainder depends on $\gamma_0$ and $C_3$. We finally 
have the existence of a constant $C>0$ (depending on $b$, $\alpha$ and $M$) such that, for every $\rho$ and $\rho'$ satisfying $d(\rho,\rho')\leq 1/2$,
$$\left|\ml{L}^u(b)(\rho)-\ml{L}^u(b)(\rho')\right|\leq C\left( d(\rho,\rho')^{\gamma_0}\log\left(\frac{1}{d(\rho,\rho')}\right)+ 
d(\rho,\rho')^{\frac{(1-\gamma_0)U_{-}}{U_+}}\right),$$
which concludes the proof of the lemma.
\end{proof}

\subsubsection{End of the proof} \label{sss:step4}

 We are now in position to give the proof of proposition~\ref{p:smooth-curve}. Thanks to lemma~\ref{l:holder}, we know that the vector field $Z_{\eps}(\rho)$ defined by~\eqref{e:derivative-heps} is a $\ml{C}^{\gamma_1}$ 
vector field and our final step is then to show that $\exp Z_{\eps}$ is in fact a good approximation of $(h^{\eps})^{-1}$ for $\|\eps\|$ small enough. 
This will follow from the fact that $\frac{dh^{\eps}}{d\eps}$ is approximated by $-\frac{dZ_{\eps}}{d\eps}$ up to an error of order $\ml{O}(\|\eps\|)$. 
In fact, we define
$$v_{\eps}(\rho):=\exp_{\rho}^{-1}h^{\eps}(\rho),$$
and fix $0<\gamma_1<\gamma_c$. According to lemma~$1$ in~\cite{Mo69}, we can write that $$\exp Z_{\eps}\circ\exp v_{\eps}=\exp(Z_{\eps}+v_{\eps}+r(Z_{\eps},v_{\eps})),$$ 
where $\|r(Z_{\eps},v_{\eps})\|_{\ml{C}^0}\leq C\|Z_{\eps}\|_{\ml{C}^{\gamma_1}}\|v_{\eps}\|_{\ml{C}^0}^{\gamma_1}$ for some uniform constant $C>0$ (depending on the 
manifold, $\gamma_1$, and J). In particular, we have that, in our setting, $\|r(Z_{\eps},v_{\eps})\|_{\ml{C}^0}=\ml{O}(\|\eps\|^{1+\gamma_1})$ with the 
constant involved in the remainder which is uniform for $(x_0,\xi_0)$ in the allowed energy layers.

\begin{rema}\label{r:moser} The proof of this fact was given in the appendix of~\cite{Mo69} for the general case of vector fields on a Riemannian manifold. The only difference is that 
the proof given in this reference is for $\gamma_1=1$. Yet, the proof can be directly adapted to get the above estimate involving H\"older norms. This property will also be used in 
paragraph~\ref{ss:eq-horocycle}.
 
\end{rema}

We finally have that, for every $(x_0,\xi_0)$ in the allowed energy layers,
$$ \tilde{\alpha}_{x_0,\xi_0}(\eps)=\exp\left(Z_{\eps}+v_{\eps}+r(Z_{\eps},v_{\eps})\right)(\alpha_{x_0,\xi_0}(\eps)).$$
Combining the definition of $Z_{\eps}$ to~\eqref{e:final-RHS}, we know that $Z_{\eps}+v_{\eps}$ is close to the zero vector field up to an error of order $\ml{O}(\|\eps\|^2)$. 
As stated above, the other vector field in the exponential map is of order $\ml{O}(\|\eps\|^{1+\gamma_1})$. Thus, we can conclude that, for every $\eps\in(-b_0,b_0)^{J+1}$ and 
uniformly for $(x_0,\xi_0)$ (in the allowed energy layers),
$$d_{T_{1/2}^*M}(\tilde{\alpha}_{x_0,\xi_0}(\eps),\alpha_{x_0,\xi_0}(\eps))=\ml{O}(\|\eps\|^{1+\gamma_1}).$$
This property holds for any $\gamma_1<\gamma_c$ which concludes the proof of proposition~\ref{p:smooth-curve}.

\begin{rema}
We have been able to construct a $\ml{C}^1$ map $\eps\mapsto\tilde{\alpha}_{x_0,\xi_0}(\eps)$ which is 
$\ml{O}(\|\eps\|^{1+\gamma_1})$ to the map $\eps\mapsto\alpha_{x_0,\xi_0}(\eps)$ we will have to consider in section~\ref{s:mixing}. The fact that the error is of 
order $\ml{O}(\|\eps\|^{1+\gamma_1})$ with $\gamma_1>0$ (and not only $o(\|\eps\|)$) will be important in the upcoming dynamical argument.  
\end{rema}

\section{Equidistribution of unstable manifolds}

\label{s:mixing}

Recall that, thanks to proposition~\ref{p:semiclassical-reduction}, 
the proof of theorem~\ref{t:maintheo} is reduced to the study of the convergence of the following integrals
\begin{equation}\label{e:classical-average}
I_{x_0,\xi_0}(b_0,T_0):=\frac{1}{(2b_0)^{J+1}}\int_{(-b_0,b_0)^{J+1}}\tilde{a}\circ G_{\eps}^{T_0}(x_0,\xi_0)d\eps,
\end{equation}
when $T_0\rightarrow+\infty$, $b_0\rightarrow 0$ (simultaneously), $p_0(x_0,\xi_0)\in [E_1-\delta/2,E_2+\delta/2]$ (for some small $\delta>0$), 
and $\tilde{a}$ is the $0$-homogeneous extension to $T^*M-M$ of a smooth function $a$ on $T_{1/2}^*M$. Recall also that $G_{\eps}^t$ is the 
Hamiltonian flow associated to $p_{\eps}(x,\xi)=p_0(x,\xi)+ V(\eps,x)$ with
$$V(\eps,x)=\eps_0 V_0(x)+\eps_1V_1(x)+\ldots +\eps_JV_J(x).$$ 
Because of the statement of proposition~\ref{p:semiclassical-reduction}, we also need to prove some uniform convergence with respect to the ``parameters'' $(x_0,\xi_0)$.

Our main statement on this purely dynamical question is given by the following theorem:

\begin{theo}\label{t:dynamics} Let $M$ be a smooth compact Riemannian surface with constant curvature $K\equiv -1$. Suppose $(V_j)_{j=0,\ldots J}$ 
satisfies~\eqref{e:admissibility}. Let $0<\gamma_1<1/2$. Let $0<\delta<E_1\leq E_2<+\infty$, and let $a$ be a smooth function on $T_{1/2}^*M$.

Then, there exist $C>0$, $0<b_1<1$, $T_1>0$, and a nonincreasing function $R(T)\rightarrow 0$ as $T\rightarrow +\infty$ such that, for every 
$(x_0,\xi_0)$ satisfying $p_0(x_0,\xi_0)\in [E_1-\delta/2,E_2+\delta/2]$, for every $b_0\in(0,b_1]$, and for every $T_0\geq T_1$, one has
\begin{equation}
 \left|\frac{1}{(2b_0)^{J+1}}\int_{(-b_0,b_0)^{J+1}}\tilde{a}\circ G_{\eps}^{T_0}(x_0,\xi_0)d\eps-\int_{T_{1/2}^*M}adL\right|\leq C
\left(b_0^{1+\gamma_1}e^{T_0\sqrt{2E}}+R\left(b_0e^{T_0\sqrt{2E}}\right)+b_0\right),
\end{equation}
where $E:=p_0(x_0,\xi_0)$.
\end{theo}

\begin{rema}\label{r:convergence-int}
At first sight, it is maybe not completely obvious if we can ensure that the upper bound in this theorem goes to $0$ for a good choice of parameters 
$b_0$ and $T_0$. Yet, we note that the term $b_0^{1+\gamma_1}e^{T_0\sqrt{2E}}$ is smaller than the term 
$b_0e^{T_0\sqrt{2E}}$. In order to ensure convergence, we have to consider a scale of parameters for which the first one goes to $0$, 
and the second one to $+\infty$. This can be achieved if we choose $E_2-E_1+2\delta$ to be not too big.



\end{rema}

More specifically, we can write the following corollary:
\begin{coro}\label{c:maincoro} Let $M$ be a smooth compact Riemannian surface with constant curvature $K\equiv -1$. Suppose $(V_j)_{j=0,\ldots J}$ 
satisfies~\eqref{e:admissibility}. 

Then for every $1<c_1\leq c_2<\frac{3}{2}$, there exists $\delta_0>0$ such that, for every $a$ in $\ml{C}^{\infty}(T_{1/2}^*M)$, one has
$$\frac{1}{(2b_0)^{J+1}}\int_{(-b_0,b_0)^{J+1}}\tilde{a}\circ G_{\eps}^{t|\log(b_0)|}(x_0,\xi_0)d\eps\longrightarrow \int_{T_{1/2}^*M}adL,\ 
\text{as}\ \eps_0\rightarrow 0,$$
uniformly for $(x_0,\xi_0)$ satisfying $p_0(x_0,\xi_0)\in [(1-\delta_0)/2,(1+\delta_0)/2]$ and uniformly for $t$ in $[c_1,c_2]$.


\end{coro}

We emphasize that an important aspect of this result is that it holds uniformly\footnote{In remark~\ref{r:burger}, we also describe the rate of convergence using results from~\cite{Bu90}.} 
for $(x_0,\xi_0)$ in the allowed energy layers and $t$ in the allowed range of times.

This corollary tells us that, for a regime of times just slightly longer than $|\log(b_0)|$, we have some equidistribution property for the evolution of points close to 
$T_{1/2}^*M$ under a family of perturbed Hamiltonian flows. For times slightly shorter than $|\log(b_0)|$ (say $T_0<(1-\delta_1)|\log(b_0)|$), we could combine 
the proof of lemma~\ref{l:transform-integral-2} below and the intertwining formulas between the geodesic and the horocycle flows~\cite{Mar77} in order to show 
that $a\circ G_{\eps}^{T_0}(x_0,\xi_0)$ is in fact equal to $a\circ G_0^{T_0}(x_0,\xi_0)$ up to some small error terms for every $\eps\in(-b_0,b_0)^{J+1}$. Thus, $|\log(b_0)|$ is 
really the \emph{critical scale of times} for which the perturbation starts to play a role in this problem.

This equidistribution property holds up to times $T_0<(1+1/2-\delta_1)|\log(b_0)|$, and it would be of course interesting to understand if the same property holds 
for much longer times than $3|\log(b_0)|/2$. In particular, this would allow to improve theorems~\ref{t:appl} and~\ref{t:maintheo} but it would also be interesting from 
the purely dynamical point of view. The factor $1/2$ appearing here is related to the H\"older regularity of the conjugating homeomorphism $h_{x_0,\xi_0}^{\eps}$ in 
the strong structural stability theorem -- see lemma~\ref{l:holder}. 

\begin{rema} In fact, our upper bound on the allowed scales of times is due to the error term $b_0^{1+\gamma_1}e^{T_0\sqrt{2E}}$ in theorem~\ref{t:dynamics} 
which will essentially appear in lemma~\ref{l:transform-integral} because we will have to replace $\eps\mapsto (h_{x_0,\xi_0}^{\eps})^{-1}(x_0,\xi_0/\|\xi_0\|)$ by a 
$\ml{C}^1$ map in order to use the equidistribution results of the unstable manifold. Recall that, to our knowledge, the map $\eps\mapsto (h_{x_0,\xi_0}^{\eps})^{-1}$ in the strong structural stability theorem is a priori not of class $\ml{C}^1$ -- see remark~\ref{r:inv-not-smooth}.
\end{rema}

This section will be devoted to the proof of theorem~\ref{t:dynamics} which is the classical counterpart of theorem~\ref{t:maintheo}. The proof will be divided in two main steps. First, we use the results of section~\ref{s:struct-stab} in order to transform the problem into a problem of equidistribution of unstable manifolds. This is the purpose of lemmas~\ref{l:transform-integral} and~\ref{l:transform-integral-2}. Then, we use some standard results on unique ergodicity of the horocycle flow~\cite{Fu73, Mar77, Bu90} in order to conclude in the case of surfaces of \emph{constant} negative curvature.

We underline that all our proof is valid for a general surface of variable negative curvature except for this final step (paragraph~\ref{ss:constant}) but this can probably be improved modulo some extra work. 

\subsection{Application of the strong structural stability property}
\label{ss:eq-horocycle}

As mentionned above, our first step will be to apply the results of section~\ref{s:struct-stab} in order to transform the integral $I_{x_0,\xi_0}(b,T_0)$ into an integral 
that can be computed using equidistribution of unstable manifolds under the geodesic flow. Precisely, we have, using the notations of proposition~\ref{p:smooth-curve}:
\begin{lemm}\label{l:transform-integral} Let $0<\gamma_1<\gamma_c$. There exists $0<b_1<1$ and $T_1>0$ such that, 
for every $(x_0,\xi_0)$ such that $E=p_0(x_0,\xi_0)\in[E_1-\delta/,E_2+\delta/2]$, for every $b_0\in(0,b_1]$, and for every $T_0\geq T_1$, one has
 $$\left|I_{x_0,\xi_0}(b_0,T_0)-\frac{1}{(2b_0)^{J+1}}\int_{(-b_0,b_0)^{J+1}}\tilde{a}\circ G_{0}^{T_0\sqrt{2E}}(\tilde{\alpha}_{x_0,\xi_0}(\eps))d\eps\right|\leq 
Ce^{U_+T_0\sqrt{2E}}b_0^{1+\gamma_1},$$
where the constant $C>0$ does not depend on $(x_0,\xi_0)$. 
\end{lemm}

Recall that in constant curvature $\gamma_c=\frac{1}{2}$ and $U_+=1$.

\begin{proof} Recall that $\tilde{a}$ is the homogeneous extension of a smooth function defined on $T_{1/2}^*M$. 
We apply the strong structural stability property under the form of remark~\ref{r:reparametrization}; more precisely, we use equation~\eqref{e:integral-rewritten}, i.e.
$$ I_{x_0,\xi_0}(b_0,T_0)=\frac{1}{(2b_0)^{J+1}}\int_{(-b_0,b_0)^{J+1}}a \circ h_{x_0,\xi_0}^{\eps}\circ G_{0}^{\tau_{x_0,\xi_0}^{\eps}\left(T_0\sqrt{2E},x_0,\frac{\xi_0}{\|\xi_0\|}\right)}
\circ\left(h_{x_0,\xi_0}^{\eps}\right)^{-1}
\left(x_0,\frac{\xi_0}{\|\xi_0\|}\right)d\eps.$$
Using~\eqref{e:taylor-stability-equation}, we can write that
$$ I_{x_0,\xi_0}(b_0,T_0)=\frac{1}{(2b_0)^{J+1}}\int_{(-b_0,b_0)^{J+1}}a\circ G_{0}^{\tau_{x_0,\xi_0}^{\eps}\left(T_0\sqrt{2E},x_0,\frac{\xi_0}{\|\xi_0\|}\right)}
\circ\left(h_{x_0,\xi_0}^{\eps}\right)^{-1}
\left(x_0,\frac{\xi_0}{\|\xi_0\|}\right)d\eps+\ml{O}(b_0),$$
where the constant in the remainder is uniform for $(x_0,\xi_0)$ in the allowed energy layers. Using~\eqref{e:taylor-stability-equation} one more time, 
we can also observe that 
$$\tau_{x_0,\xi_0}^{\eps}\left(T_0\sqrt{2E},x_0,\frac{\xi_0}{\|\xi_0\|}\right)=T_0\sqrt{2E}+\ml{O}(b_0T_0\sqrt{E}),$$
where the constant in the remainder is still uniform for $(x_0,\xi_0)$. It allows us to write the following approximation:
$$ I_{x_0,\xi_0}(b_0,T_0)=\frac{1}{(2b_0)^{J+1}}\int_{(-b_0,b_0)^{J+1}}a\circ G_{0}^{T_0\sqrt{2E}}
\circ\left(h_{x_0,\xi_0}^{\eps}\right)^{-1}
\left(x_0,\frac{\xi_0}{\|\xi_0\|}\right)d\eps+\ml{O}(b_0 T_0\sqrt{E}).$$
As a final step, we replace the map $\displaystyle\eps\mapsto\left(h_{x_0,\xi_0}^{\eps}\right)^{-1}
\left(x_0,\frac{\xi_0}{\|\xi_0\|}\right)$ (which is a priori not $\ml{C}^1$) by the $\ml{C}^1$ map $\tilde{\alpha}_{x_0,\xi_0}$ of proposition~\ref{p:smooth-curve}, and 
we find using remark~\ref{r:bound-derivative-flow} that
 $$ I_{x_0,\xi_0}(b_0,T_0)=\frac{1}{(2b_0)^{J+1}}\int_{(-b_0,b_0)^{J+1}}a\circ G_{0}^{T_0\sqrt{2E}}(\tilde{\alpha}_{x_0,\xi_0}(\eps))d\eps
+\ml{O}(b_0 T_0\sqrt{E})+\ml{O}(e^{U_+T_0\sqrt{2E}}b_0^{1+\gamma_1}).$$
\end{proof}

Thanks to lemma~\ref{l:transform-integral}, we are now reduced to understand the ``equidistribution properties'' of the map $\eps\mapsto\tilde{\alpha}_{x_0,\xi_0}(\eps)$ 
under the unperturbed flow $G_0^t$. For that purpose, we will make use of the ergodic properties of the horocycle flow, and before that, we need to make the unstable horocycle flow appear in the above expression. Recall that, in section~\ref{s:geom-back}, we defined a certain parametrization $(H^{\tau}_u)_{\tau\in\IR}$ of the horocycle flow~\cite{Mar77}. 
Our next reduction in the computation of the integral $I_{x_0,\xi_0}(b_0,T_0)$ will be to prove:

\begin{lemm}\label{l:transform-integral-2} Let $0<\gamma_1<\gamma_c$. There exists $0<b_1<1$ and $T_1>0$ such that, 
for every $(x_0,\xi_0)$ such that $E=p_0(x_0,\xi_0)\in[E_1-\delta/,E_2+\delta/2]$, for every $b_0\in(0,b_1]$, and for every $T_0\geq T_1$, one has
 $$\left|I_{x_0,\xi_0}(b_0,T_0)
-\frac{1}{(2b_0)^{J+1}}\int_{(-b_0,b_0)^{J+1}}\tilde{a}\circ G_{0}^{T_0\sqrt{2E}}\circ H_u^{\beta_{x_0,\xi_0}^u(\eps)}\left(x_0,\frac{\xi_0}{\|\xi_0\|}\right)d\eps\right|$$
$$\leq Cb_0\left(1+e^{U_+T_0\sqrt{2E}}b_0^{\gamma_1}\right),$$
where the constant $C>0$ does not depend on $(x_0,\xi_0)$, and where, using the operator defined in~\eqref{e:admissibility-operator},
$$\beta^u_{x_0,\xi_0}(\eps):=-\frac{1}{\|\xi_0\|}\sum_{j=0}^J\eps_j\ml{L}_{x_0,\frac{\xi_0}{\|\xi_0\|}}(V_j).$$ 
\end{lemm}

Under this new form, the problem we have to understand concerns now the distribution of small pieces of unstable manifolds under the geodesic flow $G_0^t$. This will be the purpose 
of the next paragraph.

\begin{proof} According to lemma~\ref{l:transform-integral} and to proposition~\ref{p:smooth-curve}, it is sufficient to approximate the integral
$$\overline{I}_{x_0,\xi_0}(b_0,T_0):=\frac{1}{(2b_0)^{J+1}}\int_{(-b_0,b_0)^{J+1}}\tilde{a}\circ G_{0}^{T_0\sqrt{2E}}\circ \exp Z_{\eps}\left(x_0,\frac{\xi_0}{\|\xi_0\|}\right)d\eps$$
We first split the vector field $Z_{\eps}=Z_{\eps}^u+Z_{\eps}^s$ into its unstable and stable components. Thanks to remark~\ref{r:moser} and 
to lemma~\ref{l:holder}, we approximate $\exp (Z_{\eps})$ by $\exp (Z_{\eps}^s)\circ \exp (Z_{\eps}^u)$ up to an error of order 
$\ml{O}(\|\eps\|^{1+\gamma_1})$, where the constant in the remainder depends only on the potential $V_j$, on $\gamma_1$ and on the manifold $M$. In particular, thanks to remark~\ref{r:bound-derivative-flow},
$$\overline{I}_{x_0,\xi_0}(b_0,T_0)=\frac{1}{(2b_0)^{J+1}}\int_{(-b_0,b_0)^{J+1}}
\tilde{a}\circ G_{0}^{T_0\sqrt{2E}}\circ \exp Z_{\eps}^s\circ \exp Z_{\eps}^u\left(x_0,\frac{\xi_0}{\|\xi_0\|}\right)d\eps+\ml{O}(b_0^{1+\gamma_1}e^{U_+T_0\sqrt{2E}}).$$
Introduce now, for $\rho\in T_{1/2}^*M$,
$$\beta^u_{\rho}(\eps):=-\frac{1}{\|\xi_0\|}\sum_{j=0}^J\eps_j\ml{L}_{\rho}(V_j),\ \text{and}\ 
k^u_{\eps}(\rho):= H_u^{\beta^u_{\rho}(\eps)}(\rho).$$
One can verify that that the maps $\eps\mapsto k^u_{\eps}$ and $\eps\mapsto \exp Z^u_{\eps}$ are of class $\ml{C}^1$ when they are considered as applications 
from $(-b_0,b_0)^{J+1}$ to $\ml{C}^0(T_{1/2}^*M,T_{1/2}^*M)$. Moreover, both maps are equal to identity in $\eps=0$ and their derivatives coincide up to an error 
of order $\ml{O}(\|\eps\|)$. In particular, we have that $\exp Z^u_{\eps}$ can be approximated by $k^u_{\eps}$ in the $\ml{C}^0$-topology up to an error of order 
$\ml{O}(\|\eps\|^2)$. The same argument can be applied to replace $\exp Z^s_{\eps}$ by a map $k^s_{\eps}$ constructed from the stable horocycle flow. Thus, we have
$$\overline{I}_{x_0,\xi_0}(b_0,T_0)=\frac{1}{(2b_0)^{J+1}}\int_{(-b_0,b_0)^{J+1}}
\tilde{a}\circ G_{0}^{T_0\sqrt{2E}}\circ k_{\eps}^s\circ k_{\eps}^u\left(x_0,\frac{\xi_0}{\|\xi_0\|}\right)d\eps+\ml{O}(b_0^{1+\gamma_1}e^{U_+T_0\sqrt{2E}}).$$ 
As $k^s_{\eps}$ is a ``reparametrization'' of the stable horocycle flow, as $k^s_{\eps}$ is close to identity (up to an error of order $\ml{O}(\|\eps\|)$) and as $T_0\geq 0$, we can 
approximate the integral as follows
$$\overline{I}_{x_0,\xi_0}(b_0,T_0)=\frac{1}{(2b_0)^{J+1}}\int_{(-b_0,b_0)^{J+1}}
\tilde{a}\circ G_{0}^{T_0\sqrt{2E}}\circ k_{\eps}^u\left(x_0,\frac{\xi_0}{\|\xi_0\|}\right)d\eps+\ml{O}(b_0^{1+\gamma_1}e^{U_+T_0\sqrt{2E}})+\ml{O}(b_0),$$
which concludes the proof of the lemma.
\end{proof}

\subsection{Using equidistribution properties of the unstable manifold}\label{ss:constant}

We underline that, up to this point of the proof, we did not use the fact that $M$ is a surface of constant curvature and all the previous arguments are in fact valid 
in the general context of a surface of variable negative curvature. In order to conclude, we now make the restriction that $M$ has a constant sectional curvature 
$K\equiv -1$. This additional assumption allows us to use the unique ergodicity of the horocycle flow for the Liouville measure~\cite{Fu73, Mar77, Bu90} 
and the fact that in this setting, our parametrization of the horocycle flow coincides with the uniformly expanding parametrization~\cite{Mar77}.

Let us recall specifically the following classical result on horocycle flows that we will use:

\begin{theo} \label{t:horocycle-conv} \cite{Fu73, Mar77} Let $M$ be a smooth, compact, Riemannian surface with constant sectional curvature $K\equiv -1$.
Let $H^s_u$ be the parametrization of the horocycle flow on $T^*_{1/2}M$ introduced in paragraph~\ref{ss:anosov}.  For any continuous function $f$ on $T_{1/2}^*M$, we have 
that
\begin{align*}
\lim_{|T| \rightarrow + \infty} \frac{1}{2T} \int_{-T}^T f \circ H^s_u(\rho) \, ds = \int_{T^*_{1/2}M} f \, dL
\end{align*}
uniformly in $\rho\in T^*_{1/2}M$ where $L$ is the Liouville measure on $T^*_{1/2}M$.
\end{theo}

In particular, one has that
$$R(T):=\sup_{|\tau|\geq T}\sup_{\rho\in T_{1/2}^*M}\left\{\left|\frac{1}{2\tau} \int_{-\tau}^{\tau} a \circ H^s_u(\rho) \, ds-\int_{T^*_{1/2}M} a \, dL\right|\right\}\longrightarrow 0,\ \text{as}\ T\rightarrow+\infty.$$

 As in the statement of theorem~\ref{t:dynamics}, we fix $0<\gamma_1<\frac{1}{2}$. We will now apply the previous theorem to the integral appearing in the conclusion of lemma~\ref{l:transform-integral-2}. Precisely, we define:
$$\tilde{I}_{x_0,\xi_0}(b_0,T_0):=\frac{1}{(2b_0)^{J+1}}\int_{(-b_0,b_0)^{J+1}}\tilde{a}\circ G_{0}^{T_0\sqrt{2E}}\circ H_u^{\beta_{x_0,\xi_0}^u(\eps)}
\left(x_0,\frac{\xi_0}{\|\xi_0\|}\right)d\eps.$$
where 
$$\beta^u_{x_0,\xi_0}(\eps):=-\frac{1}{\|\xi_0\|}\sum_{j=0}^J\eps_j\ml{L}_{x_0,\frac{\xi_0}{\|\xi_0\|}}(V_j).$$ 
This is where we will use our admissibility assumption~\eqref{e:admissibility}. This assumption implies in particular that there exists some $c_1>0$ such that, for every $(x_0,\xi_0)\in T^*M$, one can find $0\leq j_0\leq J$ such that
$$\left|\ml{L}_{x_0,\frac{\xi_0}{\|\xi_0\|}}(V_{j_0})\right|\geq c_1.$$
For a given $(x_0,\xi_0)$ in $T^*M$, we fix such a $j_0$, and we set
$$\tilde{\beta}^{u,j_0}_{x_0,\xi_0}(\eps):=-\frac{1}{\|\xi_0\|}\sum_{j=0}^J\eps_j\ml{L}_{x_0,\frac{\xi_0}{\|\xi_0\|}}(V_j)+\frac{1}{\|\xi_0\|}\eps_{j_0}\ml{L}_{x_0,\frac{\xi_0}{\|\xi_0\|}}(V_{j_0}),$$ 
which does not depend on the variable $\eps_{j_0}$. Recall now that, in the setting of constant curvature $-1$, we have the following intertwining 
relation \cite{Mar77}: $G_{0}^{T_0\sqrt{2E}} \circ  H^s \left(x_0,\frac{\xi_0}{\|\xi_0\|}\right) =  H^{s*} \circ G_{0}^{T_0\sqrt{2E}} \left(x_0,\frac{\xi_0}{\|\xi_0\|}\right)$ 
where $s^* = e^{\sqrt{2E}T_0} s$.  Performing this intertwining and applying the corresponding change of variables with respect to the parameter $\eps_{j_0}$, we obtain
$$\tilde{I}_{x_0,\xi_0}(b_0,T_0)=\frac{1}{(2b_0)^{J}}\int_{(-b_0,b_0)^{J}}\left(\frac{1}{2\tau_{x_0,\xi_0}(b_0,T_0)}
\int_{-\tau_{x_0,\xi_0}(b_0,T_0)}^{\tau_{x_0,\xi_0}(b_0,T_0)}
\tilde{a}\circ H_u^{\tau_{j_0}}\left(\rho_0(\hat{\eps}_{j_0},T_0)\right)d\tau_{j_0}\right)d\hat{\eps}_{j_0},$$
where $\hat{\eps}_{j_0}$ means that we consider all the parameters except $\eps_{j_0}$,
$$\rho_0(\hat{\eps}_{j_0},T_0):= G_{0}^{T_0\sqrt{2E}}\circ H_u^{\tilde{\beta}_{x_0,\xi_0}^{u,j_0}(\hat{\eps}_{j_0})}
\left(x_0,\frac{\xi_0}{\|\xi_0\|}\right),$$
and
$$\tau_{x_0,\xi_0}(b_0,T_0):=b_0\frac{e^{T_0\sqrt{2E}}}{\sqrt{2E}}\ml{L}_{x_0,\frac{\xi_0}{\|\xi_0\|}}(V_{j_0}).$$
Thanks to the unique ergodicity of the horocycle flow (as well as the uniformity of the remainder in $\hat{\eps}_{j_0}$), we can deduce that
$$\tilde{I}_{x_0,\xi_0}(b_0,T_0)=\int_{S^*M}\tilde{a}dL+R\left(b_0\frac{e^{T_0\sqrt{2E}}}{\sqrt{2E}}c_1\right).$$
Combining this equality to lemmas~\ref{l:transform-integral} and~\eqref{l:transform-integral-2}, we finally derive that, uniformly for $(x_0,\xi_0)$ in a small neighborhood of $T_{1/2}^*M$, one has
$$I_{x_0,\xi_0}(b_0,T_0)=\int_{S^*M}\tilde{a}dL+\tilde{R}\left(b_0e^{T_0\sqrt{2E}}\right)+\ml{O}\left(e^{T_0\sqrt{2E}}b_0^{1+\gamma_1}\right)+\ml{O}(b_0),$$
for some function $\tilde{R}(T)\rightarrow 0$ as $T\rightarrow+\infty.$ This concludes the proof of theorem~\ref{t:dynamics}.

\begin{rema}\label{r:burger}
Even if we did not use it, we underline that the results from~\cite{Bu90} (precisely theorem $2.C$) imply that $R(T)\leq C_{\alpha, a} T^{-\alpha}$, for every $0\leq \alpha<1/2$ satisfying $\alpha(\alpha-1)\geq \lambda_1$ (with $\lambda_1<0$ the first nonzero eigenvalue of $\Delta_g$ on $M$). 
In particular, this would yield that the convergence in corollary~\ref{c:maincoro} holds at a polynomial rate in $b_0$.
\end{rema}

\appendix

\section{Background on semiclassical analysis}

\label{a:pdo}
In this appendix, we review some basic facts on semiclassical analysis that can be found for instance in~\cite{Zw12}.

\subsection{General facts}

Recall that we define on $\mathbb{R}^{2d}$ the following class of symbols:
$$S^{m,k}(\mathbb{R}^{2d}):=\left\{(a_{\hbar}(x,\xi))_{\hbar\in(0,1]}\in C^{\infty}(\mathbb{R}^{2d}):|\partial^{\alpha}_x\partial^{\beta}_{\xi}a_{\hbar}|
\leq C_{\alpha,\beta}\hbar^{-k}\langle\xi\rangle^{m-|\beta|}\right\}.$$
Let $M$ be a smooth Riemannian $d$-manifold without boundary. Consider a smooth atlas $(f_l,V_l)$ of $M$, where each $f_l$ is a smooth diffeomorphism from 
$V_l\subset M$ to a bounded open set $W_l\subset\mathbb{R}^{d}$. To each $f_l$ correspond a pull back $f_l^*:C^{\infty}(W_l)\rightarrow C^{\infty}(V_l)$ and a canonical 
map $\tilde{f}_l$ from $T^*V_l$ to $T^*W_l$:
$$\tilde{f}_l:(x,\xi)\mapsto\left(f_l(x),(Df_l(x)^{-1})^T\xi\right).$$
Consider now a smooth locally finite partition of identity $(\phi_l)$ adapted to the previous atlas $(f_l,V_l)$. 
That means $\sum_l\phi_l=1$ and $\phi_l\in C^{\infty}(V_l)$. Then, any observable $a$ in $C^{\infty}(T^*M)$ can be decomposed as follows: $a=\sum_l a_l$, where 
$a_l=a\phi_l$. Each $a_l$ belongs to $C^{\infty}(T^*V_l)$ and can be pushed to a function $\tilde{a}_l=(\tilde{f}_l^{-1})^*a_l\in C^{\infty}(T^*W_l)$. 
As in~\cite{Zw12}, define the class of symbols of order $m$ and index $k$
\begin{equation}
\label{defpdo}S^{m,k}(T^{*}M):=\left\{(a_{\hbar}(x,\xi))_{\hbar\in(0,1]}\in C^{\infty}(T^*M):|\partial^{\alpha}_x\partial^{\beta}_{\xi}a_{\hbar}|\leq C_{\alpha,\beta}\hbar^{-k}\langle\xi\rangle^{m-|\beta|}\right\}.
\end{equation}
Then, for $a\in S^{m,k}(T^{*}M)$ and for each $l$, one can associate to the symbol $\tilde{a}_l\in S^{m,k}(\mathbb{R}^{2d})$ the standard Weyl quantization
$$\Op_{\hbar}^{w}(\tilde{a}_l)u(x):=
\frac{1}{(2\pi\hbar)^d}\int_{\IR^{2d}}e^{\frac{\imath}{\hbar}\langle x-y,\xi\rangle}\tilde{a}_l\left(\frac{x+y}{2},\xi;\hbar\right)u(y)dyd\xi,$$
where $u\in\mathcal{S}(\mathbb{R}^d)$, the Schwartz class. Consider now a smooth cutoff $\psi_l\in C_c^{\infty}(V_l)$ such that $\psi_l=1$ close to the support of $\phi_l$. 
A quantization of $a\in S^{m,k}(T^*M)$ is then defined in the following way~\cite{Zw12}:
\begin{equation}
\label{pdomanifold}\Op_{\hbar}(a)(u):=\sum_l \psi_l\times\left(f_l^*\Op_{\hbar}^w(\tilde{a}_l)(f_l^{-1})^*\right)\left(\psi_l\times u\right),
\end{equation}
where $u\in C^{\infty}(M)$. This quantization procedure $\Op_{\hbar}$ sends (modulo $\mathcal{O}(\hbar^{\infty})$) $S^{m,k}(T^{*}M)$ onto the space of pseudodifferential 
operators of order $m$ and of index $k$, denoted $\Psi^{m,k}(M)$~\cite{Zw12}. It can be shown that the dependence in the cutoffs $\phi_l$ and $\psi_l$ only appears at order 
$1$ in $\hbar$ (Theorem $18.1.17$ in~\cite{Ho85} or Theorem $9.10$ in~\cite{Zw12}) and the principal symbol map $\sigma_0:\Psi^{m,k}(M)\rightarrow S^{m,k}/S^{m-1,k-1}(T^{*}M)$ is 
intrinsically defined. Most of the rules (for example the composition of operators, the Egorov and Calder\'on-Vaillancourt Theorems) that hold on 
$\mathbb{R}^{2d}$ still hold in the case of $\Psi^{m,k}(M)$. Because our study concerns the behavior of quantum evolution for logarithmic times in $\hbar$, a larger class of 
symbols should be introduced as in~\cite{Zw12}, for $0\leq\overline{\nu}<1/2$,
\begin{equation}\label{symbol}
S^{m,k}_{\overline{\nu}}(T^{*}M):=\left\{(a_{\hbar})_{\hbar\in(0,1]}\in C^{\infty}(T^*M):
|\partial^{\alpha}_x\partial^{\beta}_{\xi}a_{\hbar}|\leq C_{\alpha,\beta}\hbar^{-k-\overline{\nu}|\alpha+\beta|}\langle\xi\rangle^{m-|\beta|}\right\}.
\end{equation}
Results of~\cite{Zw12} (as Calder\'on-Vaillancourt and Egorov theorems) can also be applied to this new class of symbols. 

\subsection{Egorov theorem}\label{ss:egorov}

The information in the first half of this section is largely taken from \cite{AnNo07, DyGu14}, both of which state the Egorov theorem using dynamical quantities; please 
see \cite{BaGrPa99, BoRo02} for earlier related results.

Let $e^{-\frac{it}{\hbar}\hat{P}_{0}(\hbar)}$ be the solution operator for the time dependent semiclassical Schr\"odinger equation (\ref{e:schrodinger}), which 
quantizes the flow $G_0^t$ on $T^*M$.  For a fixed energy slab 
$$T^*_{[E_1-\delta,E_2+\delta]}M:=\left\{(x,\xi)\in T^*M:p_0(x,\xi)\in[E_1-\delta,E_2+\delta]\right\},$$ 
with $0<2\delta<E_1\leq E_2<+\infty$, we consider a smooth observable $a \in C^{\infty}_c(T^*M)$ whose support is included in the energy slab $T^*_{[E_1 - \delta, E_2 + \delta]} M$.

Define the \textit{maximal expansion rate of} $G_0^t$ on the energy layer $T^*_E M$ as 
$$\Lambda^{E}_{max}:= \limsup_{|t| \rightarrow + \infty} \frac{1}{|t|} \log \sup_{\rho \in T^*_{E}M} \| d_{\rho} G_0^t \|.$$  Here, $\| d_{\rho} G_0^t \|$ is 
the operator norm of the differential $d_{\rho} G_0^t: T_{\rho} T^*M \rightarrow T_{G_0^t(\rho)} T^*M$ with respect to Sasaki metric on $T^*M$.  
By homogeneity arguments, it is possible to show that $\Lambda^E_{max}=\sqrt{2E}\Lambda^{1/2}_{max}$ and therefore deduce that the maximal expansion rate on the energy 
slab $T^*_{[E_1-\delta, E_2+\delta]}M$, defined by $\limsup_{|t| \rightarrow + \infty} \frac{1}{|t|} \log \sup_{\rho \in T^*_{[E_1 - \delta, E_2 + \delta]}M} \| d_{\rho} G_0^t \|$ is equal to $\sqrt{2E_2 + 2 \delta}\Lambda^{1/2}_{max}$.

For a fixed $T_0>0$, a fixed $\Lambda_1 > \sqrt{2E_2 + 6 \delta} \Lambda^{1/2}_{max}$ and a fixed $\overline{\nu} \in [0,1/2)$, we have that for
\begin{equation}\label{e:large-time-egorov1}
 |t| \leq T_0 +\frac{\overline{\nu}}{\Lambda_1}|\log(\hbar)|,
\end{equation}
the quantum observable $e^{\frac{it}{\hbar}\hat{P}_{0}(\hbar)} \Oph(a) e^{-\frac{it}{\hbar}\hat{P}_{0}(\hbar)} = A(t) \in \Psi^{-\infty,0}_{\overline{\nu}}(M)$ and is well-approximated in the following sense:
\begin{equation}\label{e:large-time-egorov2}
 \sigma_0(A(t)) = a \circ G_0^t + \mathcal{O}_{S^{-\infty,0}_{\nu}}(h^{1-2\overline{\nu}}) 
\end{equation}
uniformly for $t$ in the range (\ref{e:large-time-egorov1}). The subscript in  
$\mathcal{O}_{S^{-\infty,0}_{\overline{\nu}}}(\hbar^{\alpha})$ is meant to represent that the big-O estimate is taken with respect to the seminorms on
 $S^{-\infty,0}_{\overline{\nu}}(T^*M)$ with a prefactor of $\hbar^{\alpha}$.

It is worth describing the main idea of the proof
.  Given the operator equation 
$$e^{\frac{it}{\hbar}\hat{P}_{0}(\hbar)} \Oph(a) e^{-\frac{it}{\hbar}\hat{P}_{0}(\hbar)} - \Oph(a \circ G^t_0) $$ $$= \int_0^t e^{\frac{it}{\hbar}\hat{P}_{0}(\hbar)} 
\left(
\frac{i}{\hbar}\left[\hat{P}_{0}(\hbar),\Oph(a\circ G_{0}^{t-s})\right]-\Oph\left(\left\{p_{0},a\circ G_{0}^{t-s}\right\}\right)\right) e^{-\frac{it}{\hbar}\hat{P}_{0}(\hbar)}
 \, ds, $$ 
we need to show that the righthand side is a pseudo-differential operator and that its symbol lies inside the space 
$S^{-\infty,2\overline{\nu}-1}_{\overline{\nu}}(T^{*}M)$ for the prescribed time scale.  

From this formula, one can verify with the use of composition formula for pseudodifferential operators that our statement holds if we are able to control 
the $\ml{C}^k$ norms of $a\circ G_0^t$ in terms of $t$ (recall that $a$ is compactly supported in $T^*M$). Improving on earlier arguments 
from~\cite{BaGrPa99, BoRo02}, this estimate on the derivatives was done in Section 5 of \cite{AnNo07} and Appendix $C$ of~\cite{DyGu14} where it was proved that 
$$ \| a \circ G^t_{0} \|_{\mathcal{C}^k} \leq C(k) e^{k|t|\Lambda_1}  \| a \|_{\mathcal{C}^k},$$ for some positive constant $C(k)$.

In our article, we are interested in perturbations of the Schr\"odinger equation where the perturbation parameter $\eps$ lies in a small neighborhood of $0$. 
The above operator equation still holds for the perturbed problem; thus, we will be able to conclude that the Egorov theorem holds up to comparable times 
if we can control the growth of the $\ml{C}^k$ norms of $a\circ G_{\eps}^t$ in terms of similar quantities. For that purpose, we will follow the proofs of these 
references and we will start by estimating, 
uniformly for $\eps$ in a small neighborhood of $0$ and uniformly for $\rho$ in a small neighborhood of the energy slab, the quantity 
$ \|d_{\rho}G_{\eps}^t\|$ in terms of $t$, $\Lambda_{\max}^{1/2}$, $\delta$, $E_1$ and $E_2$. 


Recall that our symbol $a$ is compactly supported in $T^*_{[E_1 - \delta, E_2 + \delta]} M$. Using the definition of the maximal expansion rate of $G_0^t$, we 
know that there exists $|t_0|$ large enough, depending on the energy $E_2$ and $\delta$, 
such that for all $|t| \geq |t_0|$, $$\sup_{\rho \in T^*_{[E_1-2\delta, E_2 + 2\delta]}M} \| d_{\rho} G_0^{t} \| \leq e^{|t|\sqrt{2E_2 + 5 \delta}\Lambda_{max}^{1/2}}.$$ 

Since $\eps\mapsto G_{\eps}^{t_0}\in\ml{C}^1(T^*_{[E_1-2\delta, E_2 + 2\delta]}M)$ is $\mathcal{C}^1$ in $\eps$, it immediately follows that, there exists 
$b_0(\delta)>0$ such that, for $ \|\eps\|<b_0(\delta)$, 
$$\sup_{\rho \in T^*_{[E_1-2\delta, E_2 + 2\delta]}M} \| d_{\rho} G_{\eps}^{t_0} \| \leq \sup_{\rho \in T^*_{[E_1-2\delta, E_2 + 2\delta]}M} 
\| d_{\rho} G_0^{t_0} \| + \mathcal{O}_{t_0}(\|\eps\|)\leq e^{|t_0|\sqrt{2E_2+6\delta}\Lambda_{max}^{1/2}} .$$  
Moreover, we know that  $\Sigma^{\eps} : = \cup_{t \in \mathbb{R}} G_{\eps}^t(T^*_{[E_1-\delta, E_2 + \delta]})  \subset T^*_{(E_1-2\delta, E_2 + 2\delta)}M$ 
for all $\|\eps\| < b_1(\delta)$ (where $b_1(\delta)<b_0(\delta)$ depends on $\delta$).  
Hence, one has $$\sup_{\rho \in \Sigma^{\eps} }  \| d_{\rho} G_{\eps}^{t_0} \| \leq 
\sup_{\rho \in T^*_{[E_1-2\delta, E_2 + 2\delta]}M} \|d_{\rho} G_{\eps}^{t_0} \| \leq e^{|t_0|\sqrt{2E_2+6\delta}\Lambda_{max}^{1/2} }.$$ 
We observe that, for every integer $N$, one has 
$$d_{\rho} G_{\eps}^{Nt_0}=d_{G_{\eps}^{t_0(N-1)}\rho} G_{\eps}^{t_0}\circ \ldots \circ d_{G_{\eps}^{t_0}\rho} G_{\eps}^{t_0}\circ d_{\rho} G_{\eps}^{t_0},$$ 
and that the subset $\Sigma^{\eps}$ is invariant under the flow $G^t_{\eps}$. Then, we find that, for every $t$ in $\IR$, 
$$\sup_{\rho \in \Sigma^{\eps}}\| d_{\rho} G_{\eps}^{t} \|\leq\sup_{\rho \in T^*_{[E_1-2\delta, E + 2\delta_2]}M, 0\leq s\leq t_0} 
\| d_{\rho} G_{\eps}^{s} \|\left(\sup_{\rho \in T^*_{[E_1-2\delta, E + 2\delta_2]}M} \| d_{\rho} G_{\eps}^{t_0} \|\right)^{N},$$
where $N:=[|t|]$.This implies that there exists a constant $C(\delta)>0$ such that, for every $t$ in $\IR$ and every $\eps\in(-b_1(\delta),b_1(\delta))^{J+1}$, one has
$$\sup_{\rho \in \Sigma^{\eps}}\| d_{\rho} G_{\eps}^{t} \|\leq C(\delta)e^{|t|\sqrt{2E_2+6\delta}\Lambda_{max}^{1/2} }.$$
Thanks to this upper bound on $\|d_{\rho}G_{\eps}^{t}\|$, we can in fact show through induction that, uniformly for $\eps\in(-b_1,b_1)^{J+1}$ (where $b_1>0$ depends on 
$\delta$, $E_1$ and $E_2$) and for $t$ in $\IR$, that, for every $a$ in $\ml{C}^{\infty}_c(T^*M)$ whose support is included in $T_{[E_1-\delta, E_2+\delta]}^*M$,
$$\|\partial^{\alpha}(a\circ G_{\eps}^t)\|_{\infty}\leq \tilde{C}(\alpha, E_1,E_2,\delta)e^{|\alpha||t|\sqrt{2E_2+6\delta}\Lambda_{\text{max}}^{1/2}}
\sup_{\alpha'\in\mathbb{N}^{2d}:|\alpha'|\leq|\alpha|}\|\partial^{\alpha'}a\|_{\infty}.$$
Recall that, in the case $\eps=0$, the proof of this property was given in~\cite{AnNo07, DyGu14} by using the fact that the higher order derivatives of 
$G_0^t$ are in some sense controlled by powers of $\|d_{\rho}G_{0}^{t}\|$, which itself is controlled by  $e^{|t|\Lambda_{\max}^{1/2}\sqrt{2E_2+6\delta} }$ 
(see for instance Lemma $C.1$ in~\cite{DyGu14}). Here, $\|d_{\rho}G_{\eps}^{t}\|$ is $\ml{O}(e^{|t|\sqrt{2E_2+6\delta}\Lambda^{1/2}_{\text{max}}})$, and 
the proof of~\cite{DyGu14} can be adapted almost verbatim to our setting in proving the derivatives of order $k$ of $G_{\eps}^t$ are uniformly 
bounded by $\ml{O}(e^{|t|k\sqrt{2E_2+6\delta}\Lambda^{1/2}_{\text{max}}})$ for $\eps$ in $(-b_1,b_1)^{J+1}$.

We can now state our analogous Egorov theorem for the Hamiltonian flow $G_{\eps}^t$, which is quantized by $e^{-\frac{it}{\hbar}\hat{P}_{\eps}(\hbar)}$.  

\begin{theo} \label{t:large-time-egorov3}
Let $0<2\delta<E_1\leq E_2<+\infty$. Furthermore, fix $T_0>0$, $\Lambda_1 > \sqrt{2E_2 + 6 \delta }\Lambda^{1/2}_{max}$ and $\overline{\nu} \in [0,1/2)$.  

Then, there exists $\hbar_0>0$ and $b_1>0$  such that, for all $0<\hbar \leq \hbar_0$ and for every $a \in C^{\infty}_c(T^*M)$ whose support is included in $T^*_{[E_1-\delta, E_2 + \delta]}M$, 
we have that in the range $\eps\in(-b_1,b_1)^{J+1}$, and 
\begin{equation} \label{t:large-time-egorov4}
 |t| \leq T_0+\frac{\overline{\nu}}{\Lambda_1}|\log(\hbar)|,
\end{equation}
the quantum observable $e^{\frac{it}{\hbar}\hat{P}_{\eps}(\hbar)} \Oph(a)e^{-\frac{it}{\hbar}\hat{P}_{\eps}(\hbar)} = A_{\eps}(t) \in \Psi^{-\infty,0}_{\overline{\nu}}(M)$ 
and is well-approximated in the following sense:
\begin{equation}
 \sigma_0(A_{\eps}(t)) = a \circ G_{\eps}^t + \mathcal{O}_{S^{-\infty,0}_{\overline{\nu}}}(\hbar^{1-2\overline{\nu}}). 
\end{equation}
Moreover, all the involved semi-norms are uniformly bounded in terms of $\eps\in(-b_1,b_1)^{J+1}$ and of $t$ 
in the above range (\ref{t:large-time-egorov4}).
\end{theo}

\subsection{Positive quantization}\label{ss:antiwick}

Even if the Weyl procedure is a natural choice to quantize an observable $a$ on $\mathbb{R}^{2d}$, it is sometimes preferrable to use a quantization 
procedure $\Op_{\hbar}^+$ that satisfies the property~: $\Op_{\hbar}^+(a)\geq 0$ if $b\geq0$. This can be achieved thanks to the anti-Wick procedure $\Op_{\hbar}^{AW}$, 
see~\cite{HeMaRo87} for instance. For $a$ in $S^{0,0}_{\overline{\nu}}(\mathbb{R}^{2d})$, that coincides with a function on $\mathbb{R}^d$ outside a compact subset of 
$T^*\mathbb{R}^d=\mathbb{R}^{2d}$, one has
\begin{equation}\|\Op_{\hbar}^w(a)-\Op_{\hbar}^{AW}(a)\|_{L^2}\leq C\sum_{1\leq|\alpha|\leq D}\hbar^{\frac{|\alpha|}{2}}
\|\partial^{\alpha}a\|,
\end{equation}
where $C$ and $D$ are some positive constants that depend only on the dimension $d$.
To get a positive procedure of quantization on a manifold, one can replace the Weyl quantization by the anti-Wick one in definition~(\ref{pdomanifold}). 
This new choice of quantization (that we will denote by $\Oph^+$) is positive and it is well defined for every element $a$ in $S^{0,0}_{\overline{\nu}}(T^*M)$ of the 
form $c_0(x)+c(x,\xi)$ where $c_0$ belongs to 
$S^{0,0}_{\overline{\nu}}(T^*M)$ and $c$ belongs to\footnote{Here we mean that there exists a compact subset $K\subset T^*M$ such that $\text{supp}(a_{\hbar})\subset K$ for every $0<\hbar\leq 1$.} $\mathcal{C}^{\infty}_c(T^*M)\cap S^{0,0}_{\overline{\nu}}(T^*M)$. The main observation is that, for such symbols, one 
has
\begin{equation}\label{e:positive-quantization}
\left\|\Oph^+(a)-\Oph(a)\right\|_{L^2}\leq C'\sum_{1\leq|\alpha|\leq D'}\hbar^{\frac{|\alpha|}{2}}\|\partial^{\alpha}a\|,
\end{equation}
where $C'$ and $D'$ are some positive constants that depend only on the manifold $M$ and on the choice of coordinate charts.

\section{Manifolds of mappings}\label{a:map-manifold}


In this appendix, we briefly recall some facts on the differential structure modelled on Banach spaces that one can put on $\mathcal{C}^r(X,Y)$ (when $r\geq 0$). For a 
more detailed exposition, we refer the reader to~\cite{Ee58, Ab63, Ee66, El67, Pa68}. 


We suppose that $X$ is a compact $\ml{C}^{\infty}$ manifold and that $Y$ is a $\ml{C}^{\infty}$ Riemannian manifold of finite dimension. Let $r\geq 0$. 
Our goal is to briefly recall how the topological space $\mathcal{C}^r(X,Y)$ can be endowed with a differential structure of Banach manifold. 
Following~\cite{Ee58, Ee66} -- see also~\cite{Ab63}, we will explain how to construct coordinate charts on this space using the existence of an exponential map 
induced by the Riemannian structure. 

The construction is as follows. There exists a neighborhood $\ml{D}$ of the zero section in $TY$ and a neighborhood $\ml{W}\subset Y\times Y$ of the diagonal such that 
$\text{exp}:\ml{D}\rightarrow \ml{W}$ is a smooth diffeomorphism, where $\text{exp}$ is the exponential map induced by the Riemannian structure~\cite{GaHuLa}. 
Take now $h$ in $\mathcal{C}^r(X,Y)$ and define the following homeomorphism
$$h^*\text{exp}:h^*\ml{D}\rightarrow\ml{W}_{h},\ (x,y,\eta)\mapsto(x,\text{exp}_y( \eta)),$$
where $\ml{W}_{h}$ is an open neighborhood in $X\times Y$ and $h^*\ml{D}:=\{(x,y,\eta)\in X\times\ml{D}:h(x)=y\}$. Introduce now $U_{h}\subset\ml{C}^r(X,Y)$ the set 
of maps $g$ such that $\text{Gr}(g)\subset\ml{W}_{h}$ and the corresponding map:
$$\Omega_{h}:g\mapsto \left(x\in X\mapsto\left(h^*\text{exp}\right)^{-1}(x, g(x))\right),$$
that goes from $U_{h}$ to the Banach space $\Gamma^r(h^*TY)$ of $\ml{C}^r$ sections of the canonical vector bundle $h^*TY\rightarrow X$. Thanks to the smoothness of the 
exponential map, one can prove that the transition maps are smooth in the Fr\'echet sense~\cite{Ee58} and the atlas associated to these charts is a ``natural atlas'' on $\ml{C}^r(X,Y)$. 
We can observe that we have the following identification for the tangent space:
$$T_h\ml{C}^r(X,Y)\cong\Gamma^r(h^*TY).$$

In our proof, we take $X=Y=T_{1/2}^*M$ but we have to consider slightly smaller spaces than $\ml{C}^0(T_{1/2}^*M,T_{1/2}^*M)$. Precisely, we need to introduce
the following space:
$$\mathcal{C}_{X_0}(T_{1/2}^*M):=\left\{h\in\ml{C}^{0}(T_{1/2}^*M,T_{1/2}^*M):\ \forall \rho\in T_{1/2}^*M,
\ \left(\frac{d}{dt}h\circ G_0^t(\rho)\right)_{t=0}=D_{X_0}h\ \text{exists}\right\},$$
where $G_0^t:T_{1/2}^*M\rightarrow T_{1/2}^*M$ is the geodesic flow. It defines a topological space with the natural topology of uniform convergence of $h$ and 
$D_{X_0} h$. The above construction also allows to put a smooth differential structure on this space~\cite{dLLMM86}; in this case, the tangent space is given, for any $h$ in 
$\mathcal{C}_{X_0}(T_{1/2}^*M)$, by the identification:
$$T_h\mathcal{C}_{X_0}(T_{1/2}^*M)\cong\Gamma_{X_0}(h^*TT_{1/2}^*M):=\Gamma^0(h^*TT_{1/2}^*M)\cap\ml{C}_{X_0}(T_{1/2}^*M,h^*TT_{1/2}^*M).$$
These manifolds appear naturally when one wants to prove strong structural stability for flows using an implicit function theorem~\cite{dLLMM86}. This is 
essentially due to the fact that the conjugating homeomorphism in the strong structural stability theorem belongs to these spaces of maps. In fact, these spaces are even too big 
to ensure the invertibility condition required to apply an implicit function theorem. This problem can solved by introducing a family of transversals to the flow~\cite{dLLMM86}. For instance, 
one can consider the normal bundle $\ml{N}:=\cup_{\rho} N_{\rho}$ defined in section~\ref{s:geom-back}. Then, one can construct a small piece of submanifold locally 
modeled on the Banach space
$$\ml{V}^0_{X_0}(T_{1/2}^*M,\ml{N}):=\ml{V}^0(T_{1/2}^*M,\ml{N})\cap\ml{C}_{X_0}(T_{1/2}^*M,TT_{1/2}^*M),$$ 
where $\ml{V}^0(T_{1/2}^*M,\ml{N})$ denotes the set of continuous vector fields on $T_{1/2}^*M$ which take values in the normal bundle $\ml{N}$ over $T_{1/2}^*M$. This 
allows to define a small piece of a submanifold $\ml{M}$ in $\mathcal{C}_{X_0}(T_{1/2}^*M)$ for which the invertibility condition is satisfied~\cite{dLLMM86} (lemma $A.7$). 
We refer 
to appendix~$A$ of~\cite{dLLMM86} for more details on this construction (in this reference, they construct a submanifold adapted to any choice of transversal and not 
only adapted to the normal transversal).

In the present article, we needed to use these manifolds of mappings in order to study the properties of the implicit equation~\eqref{e:implicit-equation} induced by 
the strong structural stability property, i.e.
$$D_{X_0}h-\tau X\circ h=0_{T_{1/2}^*M}(h),$$
when the geodesic vector field $X_0$ has the Anosov property. The function $F(X,h,\tau):=D_{X_0}h-\tau X\circ h$ is defined on the space
$$\ml{V}^2(T_{1/2}^*M)\times \ml{C}_{X_0}(T_{1/2}^*M)\times \ml{C}^0(T_{1/2}^*M,\IR).$$
In paragraphs~\ref{sss:step1} and~\ref{sss:step2}, we needed to compute the partial derivatives in each direction, and we also needed to study the invertibility of 
the map
$$(D_{2} F-D_2 0, D_3F)(X_0,\text{Id}, 1):\ml{V}_{X_0}^0(T_{1/2}^*M)\times\ml{C}^0(T_{1/2}^*M,\IR)\rightarrow\ml{V}^0(T_{1/2}^*M)\oplus\ml{V}^0(T_{1/2}^*M),$$
with the conventions of remark~\ref{r:targetspace}. These questions were treated in details in appendix $A$ of~\cite{dLLMM86} in order to prove the strong 
structural stability theorem using an implicit function theorem. These properties turnout to be central in our proof (e.g. they are responsible for the admissibility conditions~\eqref{e:admissibility-operator}); thus, for the sake of completeness, we have decided to recall briefly the arguments of~\cite{dLLMM86}, at least in the context of Anosov geodesic vector fields on negatively curved \emph{surfaces}.

Concerning the partial derivatives, one has (lemma $A.5$ in~\cite{dLLMM86}):
\begin{enumerate}
 \item $D_1 F(X_0,\text{Id}, 1):L\in\ml{V}^2(T_{1/2}^*M)\mapsto (0,-L)\in\ml{V}^0(T_{1/2}^*M)\oplus\ml{V}^0(T_{1/2}^*M)$;
 \item $D_2 F(X_0,\text{Id}, 1):Z\in \ml{V}_{X_0}^0(T_{1/2}^*M)\mapsto (Z,\ml{L}_{X_0} Z)\in\ml{V}^0(T_{1/2}^*M)\oplus\ml{V}^0(T_{1/2}^*M)$, 
where $\ml{L}_{X_0}$ denotes the Lie derivative along $X_0$;
 \item $D_3 F(X_0,\text{Id}, 1):H\in \ml{C}^0(T_{1/2}^*M,\IR)\mapsto (0,-HX_0)\in\ml{V}^0(T_{1/2}^*M)\oplus\ml{V}^0(T_{1/2}^*M)$.
\end{enumerate}
Let us recall the argument to get $D_2F$ as the two other derivatives are more direct. According to lemma $A.2$ in~\cite{dLLMM86}, it 
is sufficient to take $Z$ in $\ml{V}^1(T_{1/2}^*M)$. We denote by $\varphi^s$ 
the flow of $Z$ and write that
$$D_2 F(X_0,\text{Id}, 1).Z=\frac{d}{ds}\left(d\varphi^s.X_0-X_0\circ \varphi^s\right)\rceil_{s=0}
=\frac{d}{ds}\left(\left((\varphi^s)_*.X_0-X_0\right)\circ \varphi^s\right)\rceil_{s=0}.$$
Then, we obtain
$$D_2 F(X_0,\text{Id}, 1).Z=\frac{d}{ds}\left(0_{T_{1/2}^*M}\circ \varphi^s\right)\rceil_{s=0}+\frac{d}{ds}\left(\left((\varphi^s)_*.X_0-X_0\right)\right)\rceil_{s=0},$$
which implies the result. We now have to understand the invertibility of the map:
$$(D_2F-D_20,D_3F)(X_0,\text{Id}, 1):(Z,H)\mapsto (0,\ml{L}_{X_0} Z-HX_0).$$

In~\cite{dLLMM86}, the authors proved that this map defines a linear isomorphism from $\ml{V}^0_{X_0}(T_{1/2}^*M,\ml{N})\times\ml{C}^0(T_{1/2}^*M,\IR)$ onto 
$\{0\}\oplus\ml{V}^0(T_{1/2}^*M)$ and they gave a formula for the inverse map ; please see lemma $A.7$ and its proof. 

We now briefly recall their argument 
(and refer to this reference for more details). First, we observe that the tangent space $\ml{V}_{X_0}^0(T_{1/2}^*M)$ is too big to ensure injectivity 
($(X_0,0)$ belongs to the kernel of the linear map). Hence, we need to take a submanifold $\ml{M}$ in $\mathcal{C}_{X_0}(T_{1/2}^*M)$ for which the 
tangent map will be injective. In fact, if we suppose that $(Z,H)$ belongs to $\ml{V}^0_{X_0}(T_{1/2}^*M,\ml{N})\times\ml{C}^0(T_{1/2}^*M,\IR)$, then, as the normal 
bundle $\ml{N}$ is invariant under the geodesic flow, one has that $\ml{L}_{X_0}Z$ belongs to $\ml{V}^0(T_{1/2}^*M,\ml{N})$. In 
particular, if $\ml{L}_{X_0}Z-HX_0=0$, we find that $H=0$ and $\ml{L}_{X_0}Z=0$. The Anosov assumption implies then that $Z=0$ and we get the injectivity of the linear map.

It remains now to prove surjectivity. For that purpose, we will give an exact formula for the inverse map. Consider a vector field $X$ and write its decomposition in the 
``Anosov basis'' (see paragraph~\ref{ss:anosov}):
$$X=c_0X_0+c_uX^u_{1/2}+c_sX^s_{1/2}.$$
Then, one can consider 
\begin{equation}\label{e:inverse}
(Z,H)=\left(\int_0^{+\infty}(G_0^t)_*(c_sX^s_{1/2})dt-\int_ {-\infty}^{0}(G_0^t)_*(c_uX^u_{1/2})dt,-c_0\right),
\end{equation}
which belongs to $\ml{V}^0_{X_0}(T_{1/2}^*M,\ml{N})\times\ml{C}^0(T_{1/2}^*M,\IR)$ and which satisfies $\ml{L}_{X_0}Z-HX_0=X$.

\section*{Acknowledgements} S\'ebastien Gou\"ezel pointed to one of us the existence of an analytic proof of strong structural stability: we warmly thank him for discussions related to this theorem and its different proofs. We also thank Yves Colin de Verdi\`ere, Stephan De Bièvre and Livio Flaminio for discussions and suggestions related to different aspects of the article. This article was written while SE was a resident at the Institut des Hautes \'Etudes Scientifiques. GR is partially supported by the Agence Nationale de la Recherche through the Labex CEMPI (ANR-11-LABX-0007-01) and the ANR project GeRaSic (ANR-13-BS01-0007-01).

\end{document}